\documentclass[12pt]{elsarticle}
\usepackage{amsmath,amssymb,amsthm}
\usepackage{mathrsfs}
\usepackage[bookmarks=true,
            bookmarksnumbered=true,
            bookmarksopen=true,
            colorlinks,
            pdfborder=001,
            linkcolor=black]{hyperref}
\usepackage{bm}
\usepackage{array}
\usepackage{float}
\usepackage{color}
\usepackage{multirow}
\usepackage{lineno}
\usepackage{caption}
\usepackage{stmaryrd}
\usepackage{extarrows}
\usepackage{subfigure}

\newtheorem{thm}{Theorem}[section]
\newtheorem{lem}[thm]{Lemma}
\newdefinition{rmk}{Remark}[section]
\newdefinition{definition}{Definition}
\newtheorem{example}{Example}[section]

\numberwithin{equation}{section}
\numberwithin{figure}{section}
\numberwithin{table}{section}
\newenvironment{pf}{{\noindent\textbf{ Proof}}\quad}{\hfill $\square$}

\renewcommand{\vec}[1]{\mbox{\boldmath \small $#1$}}

\renewcommand{\emph}[1]{{\color{red} #1}}

\allowdisplaybreaks

\let\al=\alpha
\let\be=\beta
\let\pa=\partial
\let\lam=\lambda

\def\Dt{\mathcal{D}_t}

\begin{document}

\begin{frontmatter}
 \title{On the explicit two-stage fourth-order accurate time discretizations}

 \author{Yuhuan Yuan}
 \ead{1548602562@qq.com}
 %\address{HEDPS, CAPT \& LMAM, School of Mathematical Sciences, Peking University, Beijing 100871, P.R. China}
 \author{Huazhong Tang\corref{cor1}}
 \ead{hztang@math.pku.edu.cn}
 \address{Center for Applied Physics and Technology, HEDPS, and LMAM, School of Mathematical Sciences, Peking University, Beijing 100871, P.R. China}
  \cortext[cor1]{Corresponding author. Fax:~+86-10-62751801.}

\begin{abstract}
This paper continues to study the  explicit two-stage fourth-order accurate time discretizations \cite{Li-Du:2016,Yuan-Tang:2020}. % for the first-order ordinary differential equations (ODE).
By introducing variable weights,
%depending on the time step-size and the dependent and independent variables,
we propose
 a class of  more general explicit one-step two-stage  time discretizations, which are different from the existing methods, such as the Euler {methods}, Runge-Kutta methods, and multistage multiderivative methods  etc.
% They are different from the existing explicit two-stage fourth-order time discretizations \cite{Li-Du:2016,Yuan-Tang:2020}
% and  the common ODE solvers (such as the Euler methods, Runge-Kutta methods, multistage multiderivative methods etc.) because of the variable weights.
%These new and more general two-stage time discretizations may be fifth-order accurate in the sense of the local truncation error under certain conditions.
We study the absolute stability, the stability interval, and the intersection between the imaginary axis and the absolute stability region.
Our results show that our two-stage time discretizations can be fourth-order accurate conditionally, % under certain conditions.
  the absolute stability region of the proposed methods with   some special choices of the variable weights
 can be larger than that of the classical explicit fourth- or fifth-order Runge-Kutta method, and the interval of absolute stability can be almost twice as much as the latter.
Several  numerical experiments are carried out to demonstrate the performance and accuracy as well as the stability
of our proposed methods.
\end{abstract}

\end{frontmatter}

\section{Introduction} \label{sec:introduction}
The  explicit two-stage fourth-order accurate time discretizations are studied in \cite{Li-Du:2016,Yuan-Tang:2020}
and successfully applied to the nonlinear hyperbolic conservation laws.
They belong to the two-derivative Runge-Kutta {methods}, see \cite{Kastlunger:1972a,Chan:2010,Tsai:2014}.
In comparison with the explicit four-stage fourth-order accurate
Runge-Kutta method, they only calls the time-consuming exact or approximate Riemann solver and the initial reconstruction with the characteristic decomposition twice  at  each time step, which is half of the former.

For the sake of  simplicity, let us consider the initial-value problem  of the first-order ordinary differential {equation} (ODE)
%For the most part, only the scalar problem
\begin{align} \label{eq:ode01}
u'(t)=L(t, u), \ t\in [0,T]; \ \ \quad u(0) = u_0,
\end{align}
where $u$ is scalar and $L(t, u)$ is linear or  nonlinear with respect to $u$.
%The extension of the present discussion to the system is more straightforward.
%If necessary, some remarks will be added.
Assume that the solution $u$ of \eqref{eq:ode01} is a sufficiently smooth function of $t$
and $L$ is also  smooth, and
give a partition of the time interval by $t_{n+1}=t_n+\tau$,
$n\in \mathbb Z^+\cup\{0\}$,
where  $\tau$ denotes the time step-size.
The Taylor series expansion of $u$ in $t$ reads
\begin{align}\nonumber
u^{n+1}=& \ \Big(u +\tau u_t +\frac{\tau^2}{2!} u_{tt}
+   \frac{\tau^3}{3!} u_{ttt}+\frac{\tau^4}{4!} u_{tttt}\Big )^n+{\mathcal O}(\tau^5)
\\  \nonumber
%=&\Big(u +\tau u_t +\frac{\tau^2}{6} u_{tt}\Big)^n
%+{  2} \frac{\tau^2}{6} \Big( \big(  u
%+ \frac{\tau}{2}  u_{t}+\frac{\tau^2}{8} u_{tt}\big)_{tt} \Big )^n+{\mathcal O}(\tau^5).
=&\ 	 \Big(u +\tau L(t,u) +\frac{\al \tau^2}{2} \Dt L(t,u)\Big)^n
\\  &\    +  \frac{(1-\al)\tau^2}{2} \Big( \big(  u
	+ \frac{\tau}{3(1-\al)}  L(t,u)+\frac{\tau^2}{12(1-\al)} \Dt L(t,u)\big)_{tt} \Big )^n+{\mathcal O}(\tau^5), \label{eq:add-decomp}
\end{align}
where $\Dt = \pa_t  + L \pa_u$ and $\alpha$ does not depend on  $t,u$.

Based on the additive decomposition \eqref{eq:add-decomp} with $\alpha=1/3$,
the explicit two-stage fourth-order time-accurate discretization \cite{Li-Du:2016}
 %the autonomous problem $u'(t)=L(u)$
 can be implemented as follows
\begin{equation} \label{eq:2stage4order}
\begin{aligned}
%u^* = u^n + \frac{\tau}{2} L(u^n) +  \frac{\tau^2}{8}  \frac{d}{d t} L(u^n), \\
%u^{n+1} =  u^n + \tau L(u^n) + \frac{ \tau^2}{6}  \frac{d}{d t} L(u^n)  + \frac{ \tau^2}{3} \frac{d}{d t} L(u^*).
u^* & = u^n + \frac{\tau}{2}  L(t^n, u^n) +  \frac{\tau^2}{8}  (\Dt L)(t^n, u^n), \\
u^{n+1} & = u^n + \tau L(t^n, u^n)
+ \frac{\tau^2}{6} \Big[(\Dt L)(t^n, u^n) + { 2} (\Dt L)(t^n+\tau/2, u^*)  \Big],
\end{aligned}
\end{equation}
which can also be found in \cite[Section 3]{Kastlunger:1972a}, \cite[Section 3.2]{Chan:2010} and \cite[Section 1]{Tsai:2014}.
%where $\Dt = \frac{\pa }{\pa t} + L \frac{\pa }{\pa u}$.
%It's easy to verify that the stability function for \eqref{eq:2stage4order} is
%\begin{equation}\label{EQ1.10}
%R(z) = 1 + z  + \frac{1}{2}  z^2  + \frac{1 }{6} z^3 + \frac{1}{24} z^4.
%\end{equation}
For a general choice of $\alpha$ that $\alpha=\alpha( \hat{\tau} )$ is a differentiable function of $\hat{\tau} = \tau^p$,  $p\geq 1 $, and  satisfies $\alpha = 1/3 + \mathcal{O}(\hat{\tau})$  and $\alpha\neq 1$, the general two-stage fourth-order time-accurate discretization  \cite{Yuan-Tang:2020}
can be given as follows
\begin{equation} \label{eq:2stage4order-former}
%\begin{cases}
%u^* = u^n + \frac{\tau}{3(1-\alpha)} L(u^n) +  \frac{\tau^2}{12(1-\alpha)}  \frac{d}{d t} L(u^n), \\
%u^{n+1} =  u^n + \tau L(u^n) + \frac{\al \tau^2}{2}  \frac{d}{d t} L(u^n)  + \frac{(1-\alpha) \tau^2}{2} \frac{d}{d t} L(u^*).
%\end{cases}
\begin{aligned}
	u^* =& u^n + \frac{\tau}{3(1-\al)}  L(t^n, u^n) +  \frac{\tau^2}{12(1-\al)}  (\Dt L)(t^n, u^n),\\
	u^{n+1} =& u^n + \tau L(t^n, u^n)
	+ \frac{\tau^2}{2} \left[ \al  (\Dt L)(t^n, u^n) + (1-\al) (\Dt L)\left( t^n+ \frac{\tau}{3(1-\al)}, u^* \right) \right],
\end{aligned}\end{equation}
which are not mentioned in the literature.
It's easy to verify that the  stability polynomials for both two-stage schemes \eqref{eq:2stage4order} and   \eqref{eq:2stage4order-former} are
\begin{equation*}%\label{EQ1.10}
\pi(\theta,z) = \theta-\left(1 + z  + \frac{1}{2}  z^2  + \frac{1 }{6} z^3 + \frac{1}{24} z^4\right),
\end{equation*}
which is the same as that of the (classical) explicit four-stage fourth-order accurate
Runge-Kutta method.  For the absolute stability  \cite{Hairer:1993,Leveque:2007}, one requires that
$$
\left|1 + z  + \frac{1}{2}  z^2  + \frac{1 }{6} z^3 + \frac{1}{24} z^4\right|\leq 1.
$$
%http://minitorn.tlu.ee/~jaagup/uk/dynsys/ds2/num/Absolute/Absolute.html
It is worth noting that there exist some examples of inequivalent definitions of the region of absolute stability of a numerical method for ODEs in the literature\footnote{
{\tt\tiny http://vmm.math.uci.edu/ODEandCM/StabiltyRegionDefinitions/StabilityRegionDefinitions.html}}.

{\em Does there   exist any explicit two-stage fourth-order accurate time discretization with a larger region of absolute stability}?
The aim of this paper is to answer this question and to propose a class of new and more general explicit one-step two-stage  time discretizations with variable weights, which   depend on the time step-size and the dependent and independent variables. It should be emphasized that those new time discretizations can have larger absolute stability regions and intervals than the classical explicit fourth- or fifth-order Runge-Kutta method, when
  the variable weights are  specially chosen.

The paper is organized as follows.
Section \ref{sec:method} proposes the general two-stage fourth-order methods.
Section \ref{sec:ASR} discusses the absolute stability of the proposed methods.
%and gives the specific choices which permit larger step-sizes than method \eqref{eq:2stage4order} (or method \eqref{eq:2stage4order-former}).
Section \ref{sec:example} conducts several numerical experiments
to demonstrate the performance
and accuracy as well as the stability  of the proposed methods.
Conclusions are given in Section \ref{sec:conclusion}.

\section{Numerical methods} \label{sec:method}
This section proposes a class of new and more general explicit one-step two-stage  time discretizations.

Instead of the additive decomposition in \eqref{eq:add-decomp}, let us consider a more general decomposition
\begin{align}%\nonumber
u^{n+1}=&	 \Big(u +\tau L(t,u) +\frac{\al \tau^2}{2} \Dt L(t,u)\Big)^n
%\\  &
 +  \frac{\beta\tau^2}{2} \Big( \big(  u
	+ \frac{\tau}{3\beta}  L(t,u)+\frac{\tau^2}{12\beta} \Dt L(t,u)\big)_{tt} \Big )^n+{\mathcal O}(\tau^5), \label{eq:add-decomp2}
\end{align}
where $\al=\al(t^n,u^n,\tau)$ and $\beta=\be(t^n,u^n,\tau)$ are two variable weights,
  depending on the time step-size and the dependent and independent variables. Based on \eqref{eq:add-decomp2},  the new and explicit two-stage time discretization can be given
as follows
\begin{equation} \label{eq:2stage4ordernew}
 	\begin{aligned}
 	u^* =& u^n + \frac{\tau}{3\beta(t^n,u^n,\tau)}  L(t^n, u^n) +  \frac{\tau^2}{12\beta(t^n,u^n,\tau)}  (\Dt L)(t^n, u^n), \\
 	u^{n+1} =& u^n + \tau L(t^n, u^n)
 	+ \frac{\tau^2}{2} \left[ \al(t^n,u^n,\tau)  (\Dt L)(t^n, u^n) + \beta(t^n,u^n,\tau) (\Dt L)\left( t^*, u^* \right) \right],
 	\end{aligned}
\end{equation}
where
$$ t^* = t^n+ \frac{\tau}{3\be(t^n,u^n,\tau)}.
$$
The following theorem gives the accuracy of the new scheme \eqref{eq:2stage4ordernew} in the sense of truncation error.
\begin{thm}
	If the variable weights $\al(t,u,\tau)$ and $\be(t,u,\tau)$ satisfy
	\begin{equation} \label{eq:condition-al-be}
	\alpha(t^n,u^n,\tau) + \beta(t^n,u^n,\tau) = 1 + {\mathcal O}(\tau^3),  \quad \beta(t^n,u^n,\tau) = \frac{2}{3}+{\mathcal O}(\tau), %\quad  \frac{d \beta}{d t} = {\mathcal O}(\tau^2),
	\end{equation}
	then the   two-stage time discretizations \eqref{eq:2stage4ordernew}  are   of fourth-order accuracy in the sense of truncation error, i.e.,
	\begin{equation*}
	u^{n+1}=\Big(u +\tau u_t +\frac{\tau^2}{2!} u_{tt}
	+   \frac{\tau^3}{3!} u_{ttt}+\frac{\tau^4}{4!} u_{tttt}\Big )^n+{\mathcal O}(\tau^5).
	\end{equation*}
\end{thm}

%\begin{rmk}
%``General" lies in that the schemes \eqref{eq:2stage4ordernew} have two parameters $\alpha,\beta$ which depend on $t,u,\tau$.
%\end{rmk}

\begin{pf}
For the sake of brevity, we   omit all superscripts $n$, write $L(t^n,u^n)$ as $L$, and
use the subscript $u$ (resp. $t$) to stand for the partial derivative with respect to $u$ (resp. $t$), for example, $L_t$ and $L_{uu}$ stand  for {$\frac{\pa L}{\pa t}(t^n,u^n)$ and $\frac{\pa^2 L}{\pa u^2}(t^n,u^n)$}, respectively, etc.
%%What's more, without misunderstanding we write u as $u^n$,  .
%Such abbreviations are also applied to the partial derivatives of $L$ so that $L_{uu}$ means $L_{uu}(t^n,u^n)$.	
%	
The Taylor  series expansion of $(\Dt L)\left( t+\frac{\tau}{3\be}, u^*\right)$ at $(t, u)$ reads
\begin{align*}
&  (\Dt L)\left( t+\frac{\tau}{3\be}, u^*\right)
=~ (\Dt L)  + \frac{\tau}{3\be} (\Dt L)_t + (u^*-u)  (\Dt L)_u \nonumber  \\
& \quad \quad	+ \frac{1}{2} \left( \frac{\tau^2}{9\be^2}  (\Dt L)_{tt} + 2 (u^*-u) \frac{\tau}{3\be}  (\Dt L)_{ut}  +  (u^*-u)^2  (\Dt L)_{uu}  \right) + \cdots.
%\label{eq:dLstar}
\end{align*}
The hypothesis \eqref{eq:condition-al-be} implies
\begin{equation*}
\frac{\tau}{3\be}= {\mathcal O}(\tau) , \quad (u^*-u) = \frac{\tau}{3\be} \left(  L +  \frac{\tau}{4}  \Dt L \right) = {\mathcal O}(\tau).
\end{equation*}
Thus, one has
%$ (\Dt L)\left( t^n+\frac{\tau}{3\be}, u^*\right)$ can be rewritten as
\begin{align*}
& (\Dt L)\left( t+\frac{\tau}{3\be}, u^*\right)
=~ (\Dt L)  + \frac{\tau}{3\be}  (\Dt L)_t  + \frac{\tau}{3\be} \left(  L +  \frac{\tau}{4}  \Dt L \right)  (\Dt L)_u \nonumber  \\
&\quad \quad +  \frac{\tau^2}{18\be^2}   (\Dt L)_{uu} L^2+  \frac{\tau^2}{9\be^2}  (\Dt L)_{ut}L +  \frac{\tau^2}{18\be^2}  (\Dt L)_{tt} +{\mathcal O}(\tau^3).
%\label{eq: dLstar-1}
\end{align*}
Substituting  it  into \eqref{eq:add-decomp2} gives
\begin{align} \label{eq:unew1}
u^{n+1} =&~ u +  \tau  L  + \frac{\tau^2}{2} (\al + \be) \Dt L   + \frac{\tau^3}{6} \left\{ (\Dt L)_t + (\Dt L)_u L\right\} \nonumber \\
&~ + \frac{\tau^4}{24} \left\{ (\Dt L)_u \cdot (\Dt L) + \frac{3}{2\be} \left[ (\Dt L)_{uu} L^2+  2 (\Dt L)_{ut}L +   (\Dt L)_{tt} \right] \right\} +{\mathcal O}(\tau^5).
\end{align}
	
On the other hand, some manipulations can give
\begin{equation*}
\begin{cases}
\Dt  L =  L_u L+ L_t,  \quad \left( \Dt L  \right)_u =  L_{uu} L + L_{ut}+ (L_u)^2,  \quad \left( \Dt L  \right)_t = L_{ut}L + L_{tt} + L_u L_t,  \\
\left( \Dt L  \right)_{uu} = L_{uuu} L + L_{uut} + 3 L_{uu}L_u,  \quad \left( \Dt L  \right)_{ut} = L_{uut}L + L_{utt} + L_{uu} L_t + 2 L_{ut} L_u,  \\
\left( \Dt L  \right)_{tt} = L_{utt} L + L_{ttt} + 2 L_{ut}L_t  + L_{tt} L_u,
\end{cases}
\end{equation*}
and
\begin{equation*}
\begin{cases}
\Dt^2 L =  \Big[ L_{uu}L^2 + 2 L_{ut}L +  L_{tt} \Big] +  \Big[ (L_u)^2 L + L_u L_t \Big],  \\ % =  \Big[ L_{uu}L^2 + 2 L_{ut}L +  L_{tt} \Big] +  \Big[ L_u \Dt L \Big],  \\
\Dt^3 L =  \Big[ L_{uuu} L^3 +  3 L_{uut} L^2 + 3 L_{utt} L + L_{ttt}   \Big] + \Big[ (L_u)^3 L + (L_u)^2 L_t \Big] \\
\quad \quad \quad + \Big[ 3L_{uu} L _u L +3L_{uu}L_t L + 3L_{ut}L_uL +3L_{ut} L_t \Big] +\Big[ L_{uu} L_u L^2+ 2L_{ut} L_u L+ L_{tt}L_u \Big].
\end{cases}
\end{equation*}
Thus, one obtains
\begin{equation*} %\label{eq:useful-001-thm}
\begin{cases}
\Dt^2 L = (\Dt L)_t + (\Dt L)_u L, \\
\Dt^3 L =(\Dt L)_{uu} L^2+  3 (\Dt L)_{ut} L +   (\Dt L)_{tt} +  (\Dt L)_u \cdot (\Dt L).
\end{cases}
\end{equation*}
Combining it with \eqref{eq:unew1} yields
\begin{align*}
& u^{n+1} = u +  \tau  L  + \frac{\tau^2}{2} \Dt L  + \frac{\tau^3}{6} \Dt^2 L +  \frac{\tau^4}{24} \Dt^3 L \\
&\quad \quad  + \frac{\tau^2}{2} (1-\al - \be) \Dt L + \frac{\tau^4}{24} \left(1-\frac{3}{2\be} \right) \cdot  \left[ (\Dt L)_{uu} L^2+  2 (\Dt L)_{ut}L +   (\Dt L)_{tt} \right] + {\mathcal O}(\tau^5).
\end{align*}
Hence, if $\al(t,u,\tau),\,\be(t,u,\tau)$ satisfy \eqref{eq:condition-al-be}, then the explicit two-stage time discretization \eqref{eq:2stage4ordernew} is fourth-order accurate.
\end{pf}

\begin{rmk}
If $ \alpha = \frac13$ and $\beta = \frac23$, then \eqref{eq:2stage4ordernew} becomes
the     two-stage fourth-order time discretizations \eqref{eq:2stage4order} proposed in \cite{Li-Du:2016}.
If $\be = 1-\alpha$ and  $\alpha=\alpha( \hat{\tau} )$ is a differentiable function of $\hat{\tau} = \tau^p, (p\geq 1)$ and  satisfies $\alpha = 1/3 + \mathcal{O}(\hat{\tau})$, $\alpha\neq 1$,
then \eqref{eq:2stage4ordernew} becomes
 \eqref{eq:2stage4order-former} studied in \cite{Yuan-Tang:2020}. Obviously, those special constant weights
 satisfy the condition \eqref{eq:condition-al-be}.
%It's easy to verify that $ \alpha = \frac13, \,\beta = \frac23$ satisfies condition \eqref{eq:condition-al-be}.
%That's to say, the two-stage fourth-order time discretization \eqref{eq:2stage4order} in is one degenerate case of \eqref{eq:2stage4ordernew}.
%
%What's more,  also satisfies condition \eqref{eq:condition-al-be}, where
%Consequently, the proposed methods \eqref{eq:2stage4ordernew} can degenerate to the two-stage fourth-order time discretizations \eqref{eq:2stage4order-former} in
\end{rmk}

\section{Absolute stability analysis} \label{sec:ASR}
This section discusses the absolute stability of the general two-stage fourth-order time discretizations \eqref{eq:2stage4ordernew}, and gives some good choices  of the variable weights $\al$ and $\be$. Under the  hypothesis \eqref{eq:condition-al-be}, our attention will be paid to the case of that
\begin{equation}  \label{eq:case1}
\alpha + \beta = 1 + \frac{C}{60} \big(\tau L_u(t^n,u^n) \big)^3,
\end{equation}
where $C$ is constant. % and independent on $\tau, u, t$.

Consider the model problem
\begin{equation} \label{eq:model-eq}
u'(t) =\lambda u(t), \quad  u(0)=u_0,
\end{equation}
with $\operatorname{Re}(\lambda)\leq 0$.
Applying the general two-stage fourth-order methods \eqref{eq:2stage4ordernew} to the model problem \eqref{eq:model-eq} with $L(t,u)=\lambda u(t)$ gives
\begin{align*}
u^{n+1} & =u^{n} + z u^{n} + \frac{\alpha + \beta}{2} z^{2}   u^{n} + \frac{1}{6} z^3 u^{n} + \frac{1}{24} z^4 u^{n},
\end{align*}
where $z := \tau \lambda$.
Combining it with \eqref{eq:case1} gives
  the (absolute) stability region
  {$$R_A(C):=\{z\in\mathbb C: |f(z,  C)|\leq 1,~ \operatorname{Re}(z)\leq0
  \},$$}
  and
the stability interval
$$I(C) := \{ z\in \mathbb R: -1 \leq f(z,C) \leq 1 ,\,z \leq 0\},$$
where the increment function (or stability function) is defined by
  %https://www.win.tue.nl/casa/meetings/seminar/previous/_abstract020220_files/talk.pdf
\begin{equation} \label{eq:f-z-C}
f(z,  C) = 1 + z + \frac{1}{2} z^2  + \frac{1 }{6} z^3  + \frac{1}{24} z^4 + \frac{C}{120} z^5.
\end{equation}
%stability function, see http://www.staff.science.uu.nl/~frank011/Classes/numwisk/ch10.pdf
It is seen that the absolute stability region $R_A(C)$ of \eqref{eq:2stage4ordernew} is the same as that of the classical explicit   fourth-   and   fifth-order Runge-Kutta methods   when $C = 0$ and 1, respectively, and for  the model problem \eqref{eq:model-eq}, the   two-stage fourth-order time discretizations  \eqref{eq:2stage4ordernew} with  \eqref{eq:condition-al-be} and \eqref{eq:case1}  is fifth-order
accurate in the sense of truncation error if $C = 1$.

%It shall be noted that the general two-stage fourth-order time discretizations mentioned later refer to the methods \eqref{eq:2stage4ordernew} with  \eqref{eq:condition-al-be} and \eqref{eq:case1}.

\begin{figure}[htbp]
	\centering
	\includegraphics[width=0.7\textwidth]{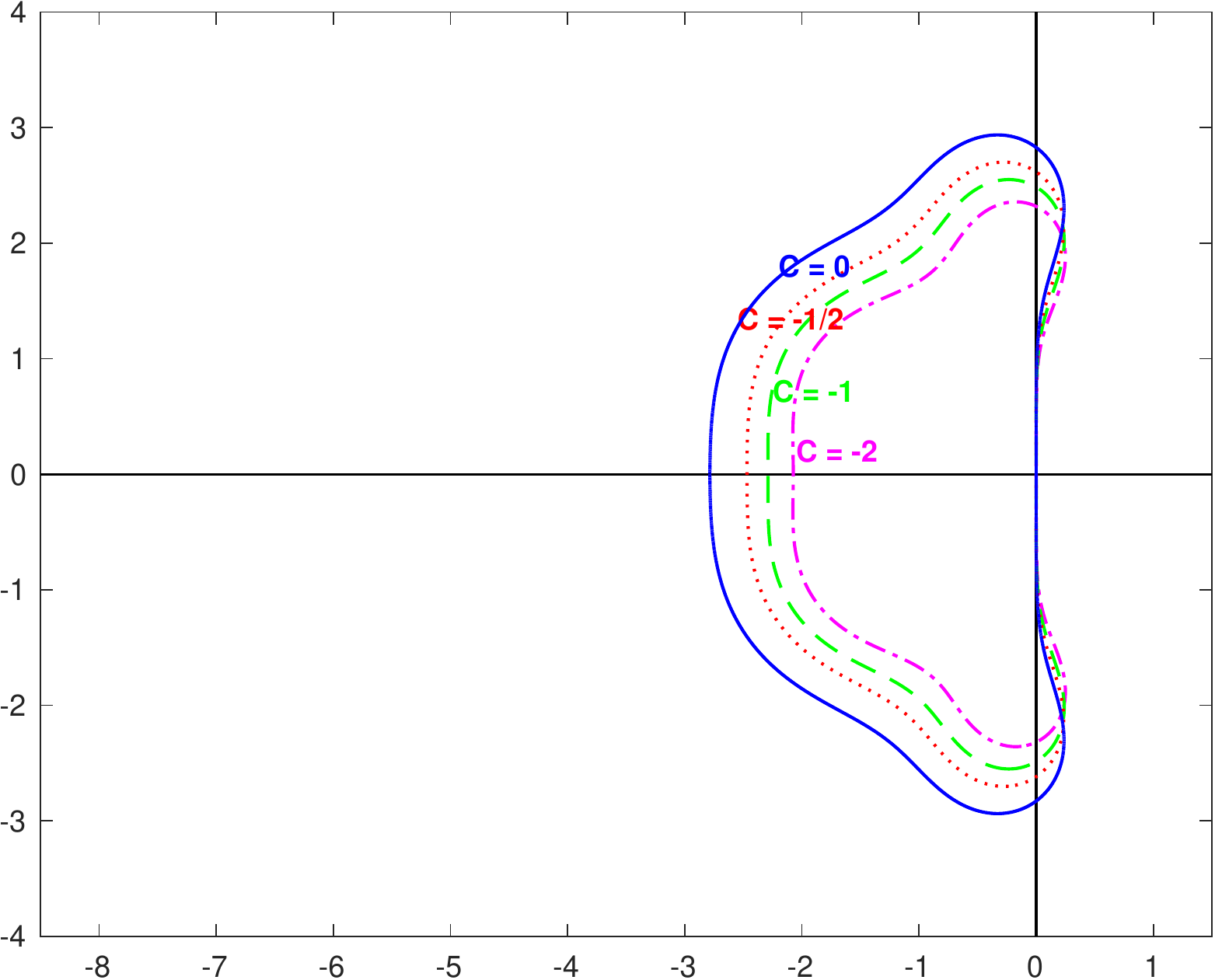}
	\caption{Curves of $|f(z,  C)|=1$ with $C = -2, -1, -\frac{1}{2}, 0$.}
	\label{fig:ASRl1}
\end{figure}

\begin{figure}[htbp]
	\centering
	\includegraphics[width=0.7\textwidth]{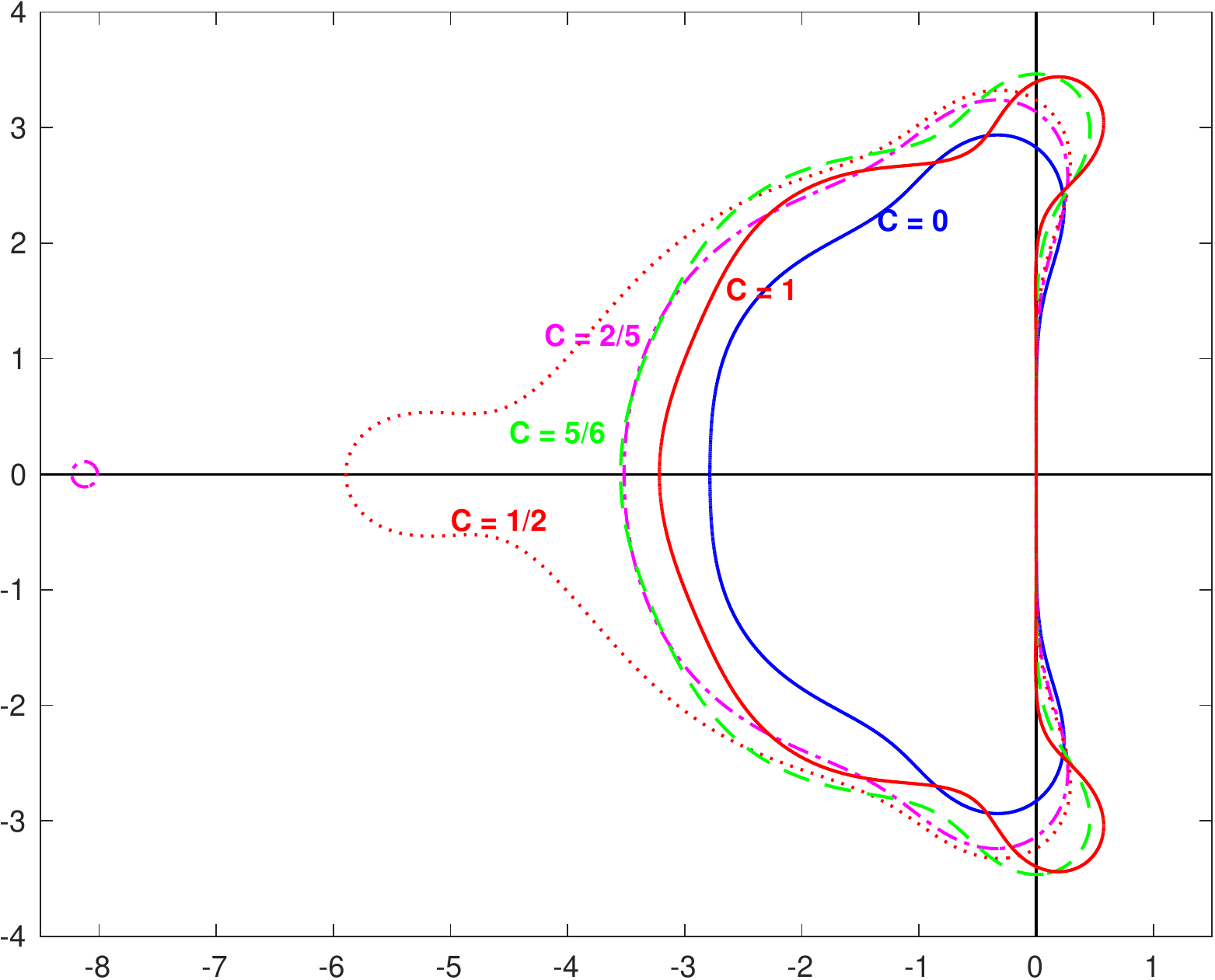}
	\caption{Curves of $|f(z,  C)|=1$ with $C = 0,\frac{2}{5}, \frac{1}{2}, \frac{5}{6}, 1$.}
	\label{fig:ASRl2}
\end{figure}

\begin{figure}[htbp]
	\centering
	\includegraphics[width=0.7\textwidth]{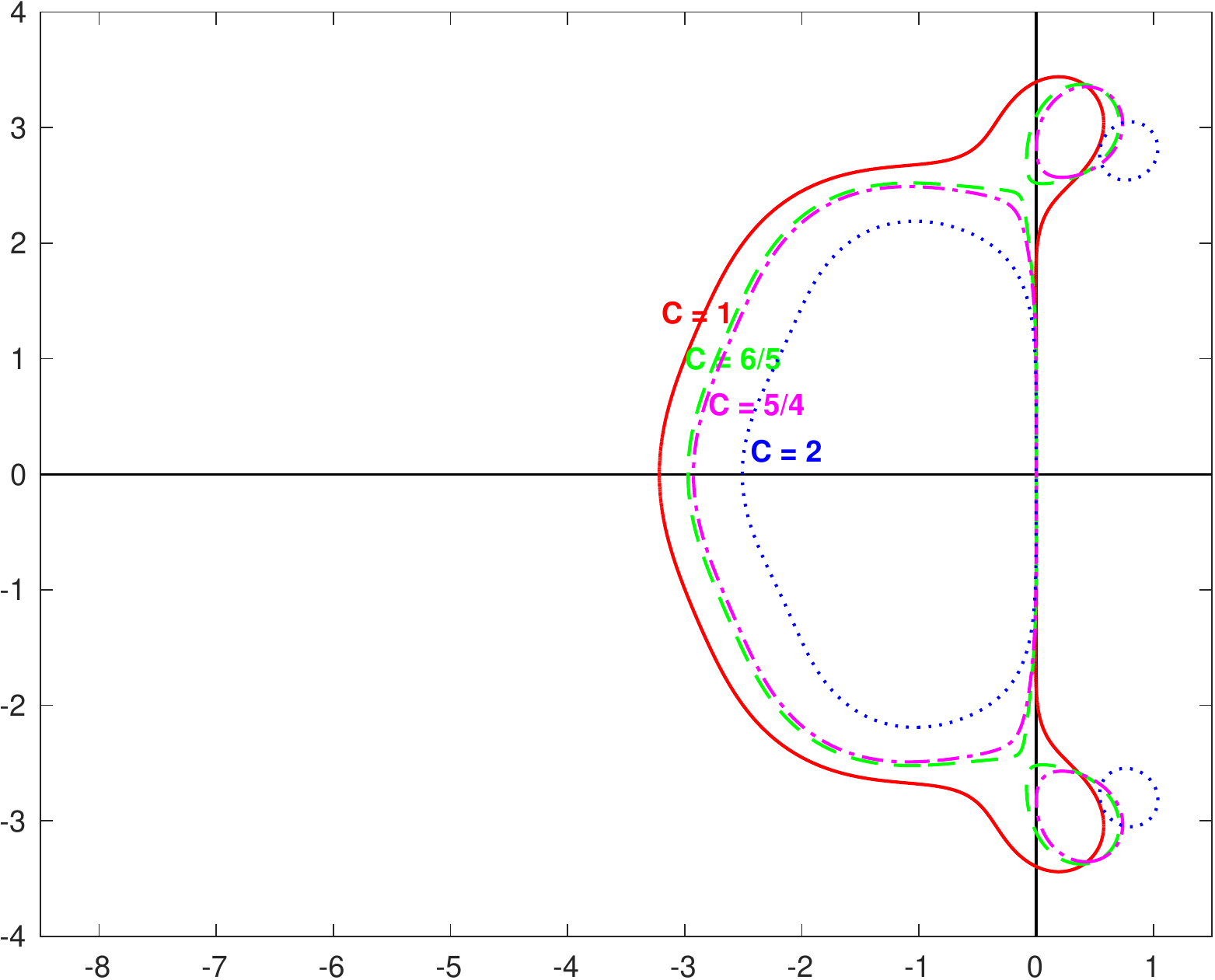}
	\caption{Curves of $|f(z,  C)|=1$ with $C = 1, \frac{6}{5}, \frac{5}{4}, 2 $.}
	\label{fig:ASRl3}
\end{figure}

Figures \ref{fig:ASRl1}-\ref{fig:ASRl3} plot the sets of complex numbers $z$ such that $|f(z,  C)|=1$, which are also showing
the  loci of the boundary of the absolute stability regions $R_A$ of the general two-stage fourth-order time discretizations \eqref{eq:2stage4ordernew} with different $C$.  %the boundary of the absolute stability regions  will .
The results show that
\begin{align*}
& \begin{cases}
R_A (-2) \subsetneq R_A(-1) \subsetneq R_A(-\frac{1}{2}) \subsetneq R_A(0),  \\
I(-2) \subsetneq I(-1) \subsetneq I(-\frac{1}{2}) \subsetneq I(0),
\end{cases}	 %\label{eq:ASR-inf-0}
\\
& \begin{cases}
R_A(0) \subsetneq R_A(\frac{2}{5}),  \  R_A(0) \subsetneq R_A(\frac{1}{2}),  \  R_A(0) \subsetneq R_A(\frac{5}{6}), \\
%I(0) \subsetneq I(1) \subsetneq I(\frac{5}{6}) \subsetneq I(\frac{1}{2}),
{ I(0) \subsetneq I(1), ~ I(1) \subsetneq I(\frac{2}{5}), ~ I(1) \subsetneq I(\frac{5}{6}) \subsetneq I(\frac{1}{2}),}
%I(0) \subsetneq I(1) \subsetneq I(\frac{1}{3}) \subsetneq I(\frac{5}{6}) \subsetneq I(\frac{1}{2}),
\end{cases} %\label{eq:ASR-inf-1}
\\
& \begin{cases}
R_A(1) \supsetneq R_A(\frac{6}{5}) \supsetneq R_A(\frac{5}{4}) \supsetneq R_A(2),  \\
I(1) \supsetneq I(\frac{6}{5}) \supsetneq I(\frac{5}{4}) \supsetneq I(2).
\end{cases} %\label{eq:ASR-1-inf}
\end{align*}
%where $R_A(C)$ denotes the absolute stability region in the left half plane,
%and $I(C)$ denotes the interval of absolute stability.
%http://minitorn.tlu.ee/~jaagup/uk/dynsys/ds2/num/Absolute/Absolute.html

\begin{rmk}
The   two-stage fourth-order time discretizations  \eqref{eq:2stage4ordernew} may be easily extended to the following system
\begin{align*} %\label{eq:ode01}
\vec u'(t)=\vec L(t, \vec u), \ t\in [0,T], \ \vec u\in \mathbb R^m,
\end{align*}
subject to $\vec u(0) = \vec u_0$, by choosing $\beta = \frac{2}{3}$ and $\vec \alpha = \frac{1}{3} \vec I_m + \frac{C  \tau^3}{60} (\nabla_{\vec u} \vec L)^3$, where $\vec I_m$ is an identity matrix
of $m\times m$.
%If $\vec u$ is a vector, then $\frac{C \tau^3}{60} (\nabla_{\vec u} \vec L )^3$ is a matrix.
%On the other hand, the value $\Dt \vec L(t^n+\frac{\tau}{3\beta}, \vec u^*)$ in our methods \eqref{eq:2stage4ordernew} hint that $\beta$ is a scalar.
%Together with \eqref{eq:case1}, one has that $\vec \alpha$ is a matrix.
%In this way,
%we can simply let
\end{rmk}

%As the stability function \eqref{eq:f-z-C} is a fourth or fifth degree polynomial in $\mathbb{C}$, it's difficult to make analysis directly on the absolute stability region.
%Hence, in the following two special cases are discussed.
%One is  the case $z \leq 0$,
%the other is the case $\operatorname{Re}(z) = 0 $.

\subsection{Interval of absolute stability}\label{subsec:Re}
%\begin{rmk}\label{rmk:Re}
This subsection discusses the interval of the absolute stability of the two-stage fourth-order time discretizations \eqref{eq:2stage4ordernew} for the case of $z=\tau\lambda \leq 0$ theoretically.
Using the definition of $f(z,C)$ and its first-order partial derivative
\begin{align}
f_z(z,  C) = 1+z+\frac12 z^2 + \frac16 z^3 + \frac{C}{24} z^4,
\label{eq:df}
\end{align}
defines
\begin{equation*}
g(z) := f(z,C) -  \frac{z}{5} f_z(z,C)= 1 + \frac{4}{5} z + \frac{3}{10} z^2 + \frac{1}{15} z^3 + \frac{1}{120} z^4.
\end{equation*}

\begin{lem}\label{lem:001}\rm
The function $g(z)$ satisfies
\begin{equation*}
g(z) > 0, \quad \mbox{for all } z\leq  0.
\end{equation*}
\end{lem}

\begin{pf}	
 By using the definition of $g(z)$,
	the derivatives of $g(z)$ are easily given as
	\begin{align*}
	g_z(z) = \frac{4}{5} + \frac{3}{5} z + \frac{1}{5} z^2 + \frac{1}{30} z^3, \
	g_{zz}(z) = \frac{3}{5}  + \frac{2}{5} z + \frac{1}{10} z^2= \frac{(z+2)^2}{10}+\frac{1}{5}.
	\end{align*}
	%which implies that $g_{zz}(z) > 0$  for all $z \in \mathbb{R}$.
	%Combining $g_z(-3) = -\frac{1}{10}< 0, g_z(-2) = \frac{2}{15} > 0$ and $g_{zz}(z) > 0$,
Because $g_z(-\infty) = -\infty< 0$, $g_z(0) = \frac{4}{5} > 0$, and $g_{zz}(z) > 0$,  $g_{z}(z)$ has a unique negative root, denoted by $z_{g_z}^*$, which is the minimum point of $g(z)$ in $(-\infty,0)$, that is, $g(z) \geq g(z_{g_z}^*)$ for all $z\leq 0$.
	Since
	\begin{equation*}
	g(z) - \left(\frac{z}{4} + \frac{1}{2} \right) g_z(z) =  \frac{3}{5}  + \frac{3}{10} z + \frac{1}{20} z^2 = \frac{(z+3)^2}{20}+\frac{3}{20} > 0,
	\end{equation*}
	one gets  $g(z_{g_z}^*) > 0$. Combining them completes the proof. %$g(z)>0$ for all $z\leq 0$.
\end{pf}

Using Lemma \ref{lem:001} yields the following conclusion.

\begin{lem}\label{rmk:fmum}\rm
The local minimum and maximum of $f(z,C)$ in $(-\infty,0)$ are positive.
\end{lem}
\begin{pf}
If using $z^*_{f_z}\in(-\infty,0)$ to denote  the negative root of $f_z(z,C)$,  then one has
\begin{equation*}
f(z^*_{f_z},C) = f(z^*_{f_z},C) -  \frac{z^*_{f_z}}{5} f_z(z^*_{f_z},C) = g(z^*_{f_z}) > 0.
\end{equation*}
The proof is completed.
\end{pf}
 \begin{figure}[htbp]
	\centering
	\includegraphics[width=0.7\textwidth]{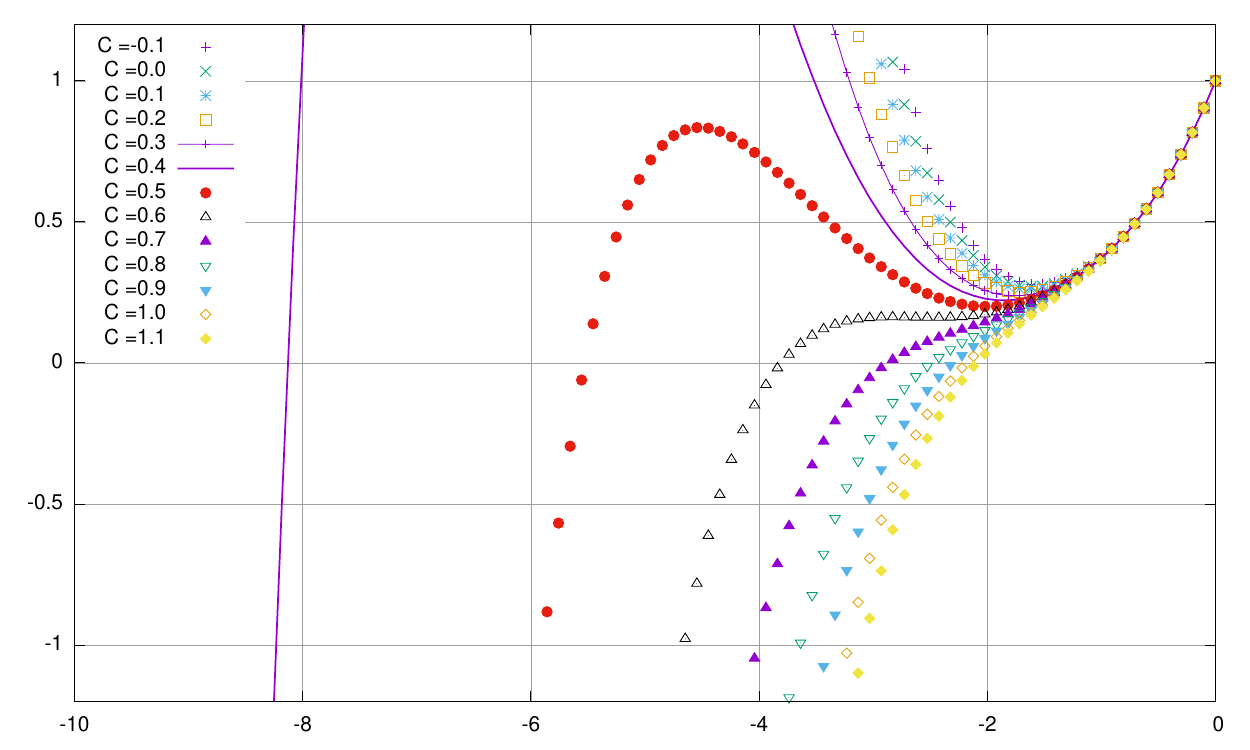}
	\caption{The profiles of $f(z,C)$.}
	\label{fig:002}
\end{figure}

In the following, we discuss the interval of the absolute stability  with the help of Lemma \ref{rmk:fmum}.
Figure \ref{fig:002} shows the profiles of $f(z,C)$ with several different $C$ in $z\in(-\infty,0)$, which can help us understand the discussion.
%\begin{description}

{\bf\tt Case 1: {$C \in (-\infty,0]$}}.
From \eqref{eq:df}, one has \begin{align}
&f_{zz}(z,  C) = 1+z+\frac12 z^2 + \frac{C}{6} z^3, \label{eq:d2f}
%&f_{zzz}(z,  C) = 1+z+ \frac{C}{2} z^2.  \label{eq:d3f}
\end{align}
thus it holds that
\begin{equation*}
f_{zz}(z,  C) \geq  1+z+\frac12 z^2  =\frac{(z+1)^2}{2}+\frac12 > 0.
\end{equation*}
It implies that $f(z,  C)$ is strictly convex for $z < 0$.
% which  means that $f(z,  C)$ changes with $z$ similar to $f(z,  0)$ shown in Figure \ref{fig:002}.
Combining {$f_{zz}(z,  C) > 0$} with $f_z(-\infty,  C)=-\infty$ and $f_z(0,  C)=1$
gives that $f_z(z, C)$ has a unique negative root, denoted by $z_{f_z}^*$, which is the minimum point of $f(z,C)$.
Using Lemma \ref{rmk:fmum} gives
\begin{equation*}
0 < f(z_{f_z}^*,C) <  f(0,C) = 1.
\end{equation*}
The readers are referred to Figure \ref{fig:002}. In this case,
for each $C \in (-\infty,0]$, the profile of $f(z,  C)$ is similar to that of $f(z,  0)$, and
the absolute stability interval $I(C)$ can be expressed as $[z^*(C), 0]$, where $z^*(C)$ is the negative solution of $f(z,C) = 1$.
With the help of the fact that $f_C(z,  C)=\frac{z^5}{120} < 0$ for $z < 0$, we can conclude that $z^*(C)$ is strictly monotonically {increasing} in $C\in (-\infty, 0]$.

{\bf\tt Case 2: $C \in (0,C_1)$}.  Here
\begin{align*}
& C_1 := \frac{-24- 24z_1 - 12 (z_1)^2 - 4 (z_1)^3}{(z_1)^4},
\ \
z_1 = -\frac{2(64+9\sqrt{67})^{1/3}}{3} + \frac{22}{3(64+9\sqrt{67})^{1/3}} -\frac83,
\end{align*}
satisfying
\begin{equation*}
f_z(z_1,C_1) = 0, \quad f(z_1,C_1) = 1.
\end{equation*}
Some computations  can show
\begin{align*}
&	z_1 \in (-4.689, -4.688), \quad C_1 \in (0.490,0.491).
\end{align*}
From \eqref{eq:d2f}, one has \begin{align*}
f_{zzz}(z,  C) = 1+z+ \frac{C}{2} z^2,  %\label{eq:d3f}
\end{align*}
and  $f_{zzz}(z,C) $ has two real roots, denoted by $z_{f_{zzz},1}^*(C)$ and $z_{f_{zzz},2}^*(C)$ with   $z_{f_{zzz},1}^*(C) < z_{f_{zzz},2}^*(C) < 0$. Since
\begin{equation*}
f_{zz}(z, C) - \frac{z}{3} f_{zzz}(z,C) = 1 + \frac{2}{3} z + \frac{1}{6} z^2 = \frac{(z+2)^2}{6} + \frac13 > 0,
\end{equation*}
the local minimum of $f_{zz}(z,C)$ satisfies $f_{zz}(z_{f_{zzz},2}^*(C), C) > 0$.
 Combining it with $f_{zz}(-\infty, C) = -\infty$ yields that $f_{zz}(z,  C)$ has only one root in $(-\infty,0)$, denoted by $z_{f_{zz}}^*(C)$, which implies that $f_z(z,C)$ is monotonically decreasing in $(-\infty, z_{f_{zz}}^*(C))$ and monotonically increasing in $( z_{f_{zz}}^*(C),0)$.
From Remark \ref{rmk:Cleq05} in the following, one has
\begin{equation} \label{eq:dfmin}
f_{z}(z_{f_{zz}}^*(C), C) < 0,	
\end{equation}
which means that $f_{z}(z,C)$ has two negative roots, denoted by $z_{f_{z},1}^*(C)$ and $z_{f_{z},2}^*(C)$ with $z_{f_{z},1}^*(C)<z_{f_{z},2}^*(C)$.
{It is worth noting that that in fact $z_1$ is a maximum point of $f(z,C_1)$, because of $f(z_{f_{z},2}^*(C_1),C_1) < f(0,C_1) = 1$.}

On the one hand, one has
\begin{equation*}
0 <	f(z_{f_{z},2}^*(C),C) < f(0,C) = 1.
\end{equation*}
On the other hand, together with $f_C(z,  C)= \frac{z^5}{120} < 0$ for $z < 0$, one has
\begin{equation*}
{f(z_{f_{z},1}^*(C_1),C)} >  f(z_{f_{z},1}^*(C_1),C_1) = f(z_1,C_1) = 1.
\end{equation*}
In this case, the profile of $f(z,  C)$ is similar to that of $f(z,  0.4)$ as shown in Figure \ref{fig:002}, and the stability interval $I(C)$ can be expressed as $[z^{*,1}(C), z^{*,2}(C)] \cup [z^{*,3}(C), 0]$, where $z^{*,2}(C),z^{*,3}(C)$,
$(z^{*,2}(C) < z^{*,3}(C))$ are the negative solutions of $f(z,C) = 1$ and $z^{*,1}(C)$ is the negative solution of $f(z,C) = -1$.

{\bf\tt Case 3:  $C \in [C_1, C_2)$}. Here
\begin{align*}
& C_2 := \frac{-6- 6z_2 - 3 (z_2)^2 }{(z_2)^3},
\ \
 z_2 = -(2+2\sqrt{3})^{1/3} + \frac{2}{(2+2\sqrt{3})^{1/3}} -2,
\end{align*}
satisfying
\begin{equation*}
f_z(z_2,C_2) = 0, \quad f_{zz}(z_2,C_2) = 0.
\end{equation*}
Similarly, by some computations, one can show
\begin{align*}
&	z_2 \in (-2.626, -2.625), \quad C_2 \in (0.603,0.604).
\end{align*}

\begin{itemize}
	\item If $C \in [C_1,0.5)$, then the same analysis for $C \in (0,C_1)$ can  give that $f_{z}(z,C)$ has two negative roots, denoted by $z_{f_{z},1}^*(C)$ and $z_{f_{z},2}^*(C)$ with $z_{f_{z},1}^*(C)<z_{f_{z},2}^*(C)$.
	
	\item If $C \in [0.5,C_2)$, then  $f_{zzz}(z,  C) \geq 0$ for any $z \in \mathbb{R}$. It means that $f_{zz}(z,  C)$ is monotonically increasing in $(-\infty,0)$.
	Combining it with $f_{zz}(-\infty,  C) = -\infty$ and $f_{zz}(0,  C) = 1$ gives that $f_{zz}(z,C)$ has a unique root in $(-\infty,0)$, denoted by  $z_{f_{zz}}^*(C) < 0$,   such  that $f_{z}(z,  C)$ is monotonically decreasing in $(-\infty,  z_{f_{zz}}^*(C))$ and monotonically increasing in $( z_{f_{zz}}^*(C),0)$.
	 Combining those with
	 $$f_{z}(-\infty,C) = +\infty, \quad f_{z}(0,C) = 1,$$
	 and
	 \begin{equation*}
{	f_{z}(z_{f_{zz}}^*(C),C) \leq 	f_{z}(z_{f_{zz}}^*(C_2),C) < f_{z}(z_{f_{zz}}^*(C_2),C_2) = f_{z}(z_2,C_2) = 0,}
	 \end{equation*}
	yields that  $f_{z}(z,C) $ has two negative roots  in $ (-\infty,0)$, denoted by $z_{f_{z},1}^*(C)$ and $z_{f_{z},2}^*(C)$ with  $z_{f_{z},1}^*(C) < z_{f_{z},2}^*(C)$.
\end{itemize}
Hence, $z_{f_{z},1}^*(C)$ is the local maximum point  of $f(z,C)$  and $z_{f_{z},2}^*(C)$ is the local minimum point of $f(z,C)$.
Together with {Lemma \ref{rmk:fmum}} and the definition of $C_1$, one can finally obtain
\begin{equation*}
{ 0 <	f(z_{f_{z},2}^*(C),C) <  f(z_{f_{z},1}^*(C),C) \leq  f(z_{f_{z},1}^*(C),C_1) \leq 1.}	
\end{equation*}
Therefore, in this case, the profile of $f(z,  C)$ is similar to $f(z,  0.5)$ as shown in Figure \ref{fig:002}, and
 the stability interval $I(C)$ can be expressed as $[z^*(C), 0]$, where $z^*(C)$ is the negative solution of $f(z,C) = -1$, and  $I(C)$ is strictly monotonically decreasing in $[C_1,C_2)$.

{\bf\tt Case 4: $C \in [C_2, \infty)$}.
Using the same analysis as that for $C \in [0.5,C_2)$ can   give that
$f_{zz}(z,C)$ has a unique root   in $(-\infty,0)$, denoted by $z_{f_{zz}}^*(C) $, and $f_{z}(z,  C)$ is monotonically decreasing in $(-\infty,  z_{f_{zz}}^*(C))$ and monotonically increasing in $( z_{f_{zz}}^*(C),0)$.
With the definition of $C_2$, one has
\begin{equation*}
{f_{z}(z_{f_{zz}}^*(C),C)  \geq f_{z}(z_{f_{zz}}^*(C),C_2) \geq f_{z}(z_{f_{zz}}^*(C_2),C_2) = f_{z}(z_2,C_2) = 0,}
\end{equation*}
which implies that the function $f(z,C)$ is monotonically increasing in $(-\infty,0)$.
Therefore, in this case, the profile of $f(z,  C)$ is  similar to $f(z,  1)$ shown in Figure \ref{fig:002}, and
the stability interval $I(C)$ can be expressed as $[z^*(C), 0]$, where $z^*(C)$ is the negative solution of $f(z,C) = -1$, and $I(C)$ strictly monotonically decreases in $[C_2,  +\infty)$.
%\end{description}

\begin{rmk} \label{rmk:Cleq05}
	Let us  verify the inequality \eqref{eq:dfmin}.
	It may be proved by contradiction. Assume that { $f_{z}(z_{f_{zz}}^*(C), C) \leq 0$}, which implies that the function $f(z,C)$ is monotonically increasing. Some calculations give
	\begin{equation*}
	f(-2,C) = \frac13 - \frac{4}{15} C, \quad f(-4,C) = 5 - \frac{128}{15}C.
	\end{equation*}
	Then, one has
	\begin{equation*}
	f(-4,C) - f(-2,C) = \frac{14}{3} - \frac{124}{15} C < 0,
	\end{equation*}
	%gives $C \geq \frac{35}{62}$,
which is in contradiction with $C \in (0,0.5)$.
\end{rmk}

\begin{rmk} \label{IntervalofRe}
	If using $z^*_{0}$, $z^*_{0.5}$ and $z^*_{1}$  to denote the  solutions of $f(z,0) = 1$, $f(z,0.5) = -1$ and $f(z,1) = -1$, respectively,
	then the intervals of the absolute stability of the general two-stage fourth-order time discretizations with $C = 0,\,  0.5,\,  1$ are
	\begin{align*}
	I(0) = [z^*_0,  0], \quad I(0.5) = [z^*_{0.5},  0], \quad I(1) = [z^*_1,  0],
	\end{align*}
	respectively, where
	$z^*_{0}$, $z^*_{0.5}$ and $z^*_{1}$ satisfy
	\begin{align*}
	z^*_0 \in (-2.786,-2.785), \quad z^*_{0.5} \in (-5.894,-5.893), \quad z^*_{1} \in (-3.218, -3.217).
	\end{align*}
\end{rmk}

\subsection{Intersection between   imaginary axis and   stability region}\label{subsec:Im}
This subsection discusses the intersection between the imaginary axis and the absolute stability region, denoted by $I_{im}(C)$.
Let $z = i \zeta$ with $\zeta \in \mathbb{R}$, $i^2 = -1$.  Then one has
\begin{equation*}
f(i \zeta, C) = 1-\frac{1}{2} \zeta^2 + \frac{1}{24} \zeta^4 + i \cdot \zeta \Big( 1-\frac{1}{6} \zeta^2 + \frac{C}{120} \zeta^4 \Big),
\end{equation*}
and the value $|f(i \zeta, C)|^2$ can be calculated by
\begin{align*}
|f(i \zeta, C)|^2  = &~ \left(  1-\frac{1}{2} \eta + \frac{1}{24} \eta^2  \right)^2 + \eta \left( 1-\frac{1}{6} \eta + \frac{C}{120} \eta^2 \right)^2 \\
= &~  1-\frac{\eta^3}{72} +\frac{\eta^4}{576} + \frac{C \eta^3( C\eta^2 - 40 \eta + 240 )}{14400},
\end{align*}
where $\eta = \zeta^2 \geq 0$.
The absolute stability requires that $\eta = 0$ or
\begin{align}\label{eq:im-001}
-\frac{1}{72} +\frac{\eta}{576} + \frac{C ( C\eta^2 - 40 \eta + 240 )}{14400} \leq 0,  \quad \mbox{for }  {\eta \geq 0}.
\end{align}
If defining
\begin{align*}
g(\eta, C) :=  C^2 \eta^2 + 5(5- 8 C) \eta + 40(6 C - 5),
\end{align*}
then \eqref{eq:im-001} is equivalent to  $g(\eta,  C) \leq 0$ for $\eta \geq 0$. In the following, we discuss its solution.

By some tedious manipulations, we can yield the conditions for \eqref{eq:im-001}.
\begin{itemize}
\item If $C = 0$, then the absolute stability requires   $0 \leq \eta \leq 8$, equivalently,
$\zeta \in [ -2\sqrt{2}, 2\sqrt{2} ]$.% = [-2.8284, 2.8284]$.
%then $f(i \zeta, 0) \leq 1$.

\item If $C\neq 0$, then calculate  the discriminant of the quadratic equation $g(\eta, C) = 0$  by
\begin{equation*}
\Delta = (25- 40 C)^2 - 4C^2(240 C - 200) = 5( 5-4C) \cdot \big( 48 C^2 - 60 C + 25 \big).
\end{equation*}
Because  $48 C^2 - 60 C + 25 > 0$ for all $C\in\mathbb R$,
  the sign of $\Delta$ is determined by $5-4C$.
\begin{itemize}
\item
 If $C > \frac{5}{4}$, then $\Delta < 0$ and the equation $g(\eta,  C) \leq 0$ for $\eta \geq 0$ has no real solution,
thus one has
	\begin{equation*}
	|f(i \zeta,  C)| \leq 1 \;\;\quad \mbox{if and only if} \;\; \quad  \zeta =  0.
	\end{equation*}

\item  If $C = \frac{5}{4} $, then  $\Delta = 0$ and the equation $g(\eta,  C) =0$ for $\eta \geq 0$ has a unique solution $\eta = 8$, so that %$g(\eta,  C) \leq 0, \, (\eta \geq 0)$ gives $\eta = 8$.
%	That's to say,
	\begin{equation*}
	|f(i \zeta,  C)| \leq 1, \;\;\quad \mbox{if and only if} \;\; \quad  \zeta = \pm 2\sqrt{2}, \, 0.
	\end{equation*}
%
%
%Then $C$ can be further classified with the help of  the sigh of the leading coefficient $C^2$, discriminant $\Delta$, axis of symmetry $-\frac{25-40C}{2C^2}$ and intercept $g(0,C)=\widetilde{}40(6C-5)$.
\item If $C <\frac{5}{4} $, then  $\Delta > 0$ and the equation $g(\eta,  C) =0$ for $\eta \geq 0$ has two different real solutions
\begin{equation*}
\eta_- = \frac{-(25-40C) - \sqrt{\Delta}}{2C^2}, \quad \eta_+ = \frac{-(25-40C)+\sqrt{\Delta}}{2C^2}.
\end{equation*}
According to the sign of $g(0,C)$, our discussion is divided into three cases.
\begin{itemize}
\item If $C < \frac{5}{6}$,  then $ g(0,C)< 0$, thus $\eta_+ > 0$ and $\eta_- < 0$, so that  the inequality $g(\eta,  C) \leq 0$ for $\eta \geq 0$ requires $0 \leq \eta \leq \eta_+$.
	That is to say,
	\begin{equation*}
	|f(i \zeta,  C)| \leq 1, \;\;\quad \mbox{if and only if} \;\; \quad \zeta \in [ -\sqrt{\eta_+},  \sqrt{\eta_+} ].
	\end{equation*}
	
	\item If $C = \frac{5}{6}$, then $\eta_+ > 0,\,\eta_- = 0$, and the inequality $g(\eta,  C) \leq 0$ for $\eta \geq 0$ requires gives $0 \leq \eta \leq \eta_+$.
	That is to say,
	\begin{equation}
	|f(i \zeta,  C)| \leq 1, \;\;\quad \mbox{if and only if} \;\; \quad \zeta \in [ -\sqrt{\eta_+},  \sqrt{\eta_+} ].
	\end{equation}
		
	\item If $C \in ( \frac{5}{6},  \frac{5}{4} ) $, then  $\eta_+ > 0,\,\eta_- > 0$, and thus the inequality $g(\eta,  C) \leq 0$ for $\eta \geq 0$ requires $\eta_- \leq \eta \leq \eta_+$.
	That is to say,
	\begin{equation*}
	|f(i \zeta,  C)| \leq 1, \;\;\quad \mbox{if and only if} \;\; \quad  \zeta \in [ -\sqrt{\eta_+},  -\sqrt{\eta_-} ] \cup  [ \sqrt{\eta_-},  \sqrt{\eta_+} ] \cup 0.
	\end{equation*}

\end{itemize}

	\end{itemize}
\end{itemize}
In all, the interval $I_{im}(C)$ can be summed up as follows
\begin{equation} \label{eq:I_Im}
I_{im}(C) =
\begin{cases}
[ -\sqrt{\eta_+},  \sqrt{\eta_+} ],  & \mbox{if} \;\; C < 0; \\
[ -2\sqrt{2},  2\sqrt{2} ],  & \mbox{if} \;\; C = 0; \\
[ -\sqrt{\eta_+},  \sqrt{\eta_+} ],  & \mbox{if} \;\; C \in (0, \frac{5}{6}]; \\
[ -\sqrt{\eta_+},  -\sqrt{\eta_-} ] \cup  [ \sqrt{\eta_-},  \sqrt{\eta_+} ] \cup \{0\},  & \mbox{if} \;\; C \in ( \frac{5}{6},  \frac{5}{4} ); \\
\{ -2\sqrt{2}, \, 2\sqrt{2} ,\, 0 \},   & \mbox{if} \;\; C  = \frac{5}{4};\\
\{0\},   & \mbox{if} \;\; C  \in (\frac{5}{4}, \infty).
\end{cases}
\end{equation}
When $C = 0, 0.5, 1$, the intersections $I_{im}(C)$ between the imaginary axis and the absolute stability regions  are explicitly and respectively given by
\begin{align*}
& I_{im}(0) = [-2\sqrt{2},  2\sqrt{2}], \quad I_{im}(0.5) = \left[-\sqrt{2(\sqrt{105}-5)}, \sqrt{2(\sqrt{105}-5)}\right], \\
& I_{im}(1) = \left[ - \sqrt{ \frac{15+\sqrt{65}}{2}},  -\sqrt{ \frac{15-\sqrt{65}}{2}}\right] \cup \left[ \sqrt{ \frac{15-\sqrt{65}}{2}},  \sqrt{ \frac{15+\sqrt{65}}{2}} \right] \cup 0.
\end{align*}

\begin{rmk}
For the case of   $C < \frac{5}{6}$ and $C \neq 0$, see \eqref{eq:I_Im}, one can know that $I_{im}(C)$ depends on the positive root $\eta_{+}$ of $g(\eta,C)$.
On the other hand, thanks to $g(8,  C) = 16 C (4C-5)$,  $g(8,C) > 0$ if $C<0$ and $g(8,C) < 0$ for $C \in (0,\frac56)$. Hence, $I_{im}(C)$ satisfies
\begin{equation*}
\begin{cases}
I_{im}(C) \subsetneq I_{im}(0),  & \mbox{if} \;\; C < 0; \\
I_{im}(C) \supsetneq I_{im}(0),   & \mbox{if} \;\; C \in (0, \frac{5}{6}).
\end{cases}
\end{equation*}
For $C  \in (\frac56, \frac54)$, see \eqref{eq:I_Im},  $I_{im}(C)$ depends on the roots $\eta_{\pm}$ of $g(\eta,C)$.  Because $\pa_{C} g(\eta,  C) = 2C\eta^2 -40\eta + 240 > \frac53\eta^2 -40\eta + 240 \geq 0$,  the set $\{\eta \,|\, g(\eta,  C) \leq 0, \,  \eta \geq 0\}$ decreases as  $C$ increases, and thus one can get
\begin{equation*}
\begin{cases}
I_{im}(C) \subsetneq I_{im}(1),  & \mbox{if} \;\; C \in (1, \frac54); \\
I_{im}(C) \supsetneq I_{im}(1),   & \mbox{if} \;\; C \in (\frac{5}{6}, 1).
\end{cases}
\end{equation*}
For $I_{im}(C)$,  the choice of $C$ in $[0,  1]$ is better than $C < 0$ and $C > 1$.
\end{rmk}

%
%\begin{rmk} \label{rmk:IntervalofIm}
%
%\end{rmk}

\section{Numerical results} \label{sec:example}
This section conducts several numerical experiments to demonstrate the performance and the above analyses of the general two-stage fourth-order time discretizations
\eqref{eq:2stage4ordernew} with \eqref{eq:case1},  in comparison with the following %classical explicit
four-stage fourth-order  Runge-Kutta method (abbreviated by {\tt RK4})
\begin{equation*} %\label{eq:4th-RK}
\begin{cases}
u^{(1)} =u^{n}+\frac{1}{2} \tau  L(t^n,u^n), \\
u^{(2)} =u^{n}+\frac{1}{2} \tau  L(t^{n}+\frac{1}{2} \tau, u^{(1)}), \\
u^{(3)} =u^{n}+ \tau  L(t^{n}+ \frac{1}{2}\tau, u^{(2)}), \\
u^{n+1} =\frac{1}{3}\Big[ u^{(1)}+2 u^{(2)}+ u^{(3)}-u^{n}+\frac{1}{2} \tau  L(t^{n}+\tau, u^{(3)}) \Big].
\end{cases}
\end{equation*}
%Since the  analyses in Section \ref{subsec:Re} and Section \ref{subsec:Im}
The following will only show the numerical results obtained with $C = 0, 0.5 , 1$, which are better than the choice of $C<0$ or $C>1$ as shown in Section \ref{sec:ASR}. %Unless otherwise stated
The diagrams will be drawn with symbols ``$\circ$'', ``$+$'' and ``$\vartriangle$'', and ``$\square$'' for
the two-stage fourth-order time discretizations
\eqref{eq:2stage4ordernew} with \eqref{eq:case1}
 and $C = 0, 0.5, 1$, and {\tt RK4}, respectively.

\subsection{Scalar case} \label{subsec:scalar}
This subsection will solve several first-order ordinary differential equations by using the two-stage high-order methods \eqref{eq:2stage4ordernew} with  $\alpha = \frac{1}{3}+\frac{C}{60} (L_u \tau)^3,~\beta = \frac{2}{3}$ or $\alpha = \frac{1}{3}, ~\beta = \frac{2}{3} +\frac{C}{60} (L_u \tau)^3$. The constant $C$ is taken as $C = 0, 0.5, 1$, respectively, and
 %Unless otherwise stated,
 the relative error
 \begin{equation*}
 err(u) = \frac{|u(T) - u_{\tau}(T)|}{|u(T)|},%|u(T) - u_{\tau}(T)|
 \end{equation*}
is estimated, where $u(T)$ and $u_{\tau}(T)$  are the exact and  numerical solutions at $t = T$, respectively.

\begin{example} \label{ex001new}
	Consider the initial value problem
	\begin{equation} \label{eq:ex001new}
	u'(t) = -u,~ t\geq 0; \quad u(0) = 1,
	\end{equation}
	whose exact solution is $u(t) = \exp(-t)$.
	%Solve this problem up to time $t = 4$.
	
	Table \ref{table:acc-ex001new-largetau} gives the relative errors and convergence rates at $t=4$ obtained by the method \eqref{eq:2stage4ordernew}  with $\alpha = \frac{1}{3}+\frac{C}{60} (L_u \tau)^3, \beta = \frac{2}{3}$, where the reference step-sizes $\tau_0$ are given by the   intervals of the absolute stability in Remark \ref{IntervalofRe}. The result clearly shows that for the model problem \eqref{eq:ex001new} the proposed method with $C = 0$ or $0.5$ is fourth-order accurate, and  \eqref{eq:2stage4ordernew}  with $C = 1$ is fifth-order accurate.	It is consistent with the previous %theoretic
theoretical analysis.
In this test, the numerical results obtained by   $\alpha = \frac{1}{3},\beta = \frac{2}{3}+\frac{C}{60} (L_u \tau)^3$ %and {\tt RK4}
are the same as those in Table \ref{table:acc-ex001new-largetau}, so  they
are not presented here to avoid  repetition.
	%This is because for the model equation the numerical results obtained by our methods with $\alpha = \frac{1}{3},~\beta = \frac{2}{3}+\frac{C}{60} (L_u \tau)^3$ are the same to that with $\alpha = \frac{1}{3}+\frac{C}{60} (L_u \tau)^3,\beta = \frac{2}{3}$ and that obtained by {\tt RK4}.
%	And the numerical results obtained by {\tt RK4} are the same to that obtained by our methods with $C = 0$.	

\end{example}

\begin{table}[htpb]
	\centering
	\caption{ Example \ref{ex001new}: The relative errors and convergence rates of the general two-stage fourth-order time discretization with $\alpha = \frac{1}{3}+\frac{C}{60} (L_u \tau)^3,\,\beta = \frac{2}{3}$ and  reference step-size $\tau_0$.}
	\label{table:acc-ex001new-largetau}
	\begin{tabular}{|c|cc|cc|cc|}
		\hline%\toprule
		\multirow{2}{*}{$\tau$}  &      \multicolumn{2}{c|}{$C=0$, $\tau_0 = 2.7$} & \multicolumn{2}{c|}{$C=0.5$, $\tau_0 = 5.8$}& \multicolumn{2}{c|}{$C=1$, $\tau_0 = 3.2$}\\
		\cline{2-7}
		&   error & order   & error & order&   error & order\\
		\hline%\midrule
		
		$\tau_0 $ &   1.3291e+01 & - &  3.9039e+01 & - &   2.4742e+01 & - \\
		$\tau_0/2 $  &    3.6366e-01  &  5.1917  &  5.1269e+00  &  2.9287 &   1.7886e-01  &  7.1120 \\
		$\tau_0/4 $  &   1.1691e-02  &  4.9591 &   1.5732e-01  &  5.0263 &     3.6257e-03  &  5.6244 \\
		$\tau_0/8 $  &   5.5332e-04  &  4.4011 &   6.7895e-03  &  4.5343 &    8.0248e-05  &  5.4976 \\
		$\tau_0/16 $  &   3.0414e-05  &  4.1853 &    3.6496e-04  &  4.2175 &   2.1109e-06  &  5.2486 \\
		$\tau_0/32 $  &   1.7974e-06  &  4.0807 &    2.0228e-05  &  4.1733 &   6.0532e-08  &  5.1240 \\
		\hline%\bottomrule
	\end{tabular}
\end{table}

\begin{figure}[htpb]
	\centering
	\includegraphics[width=0.49\textwidth]{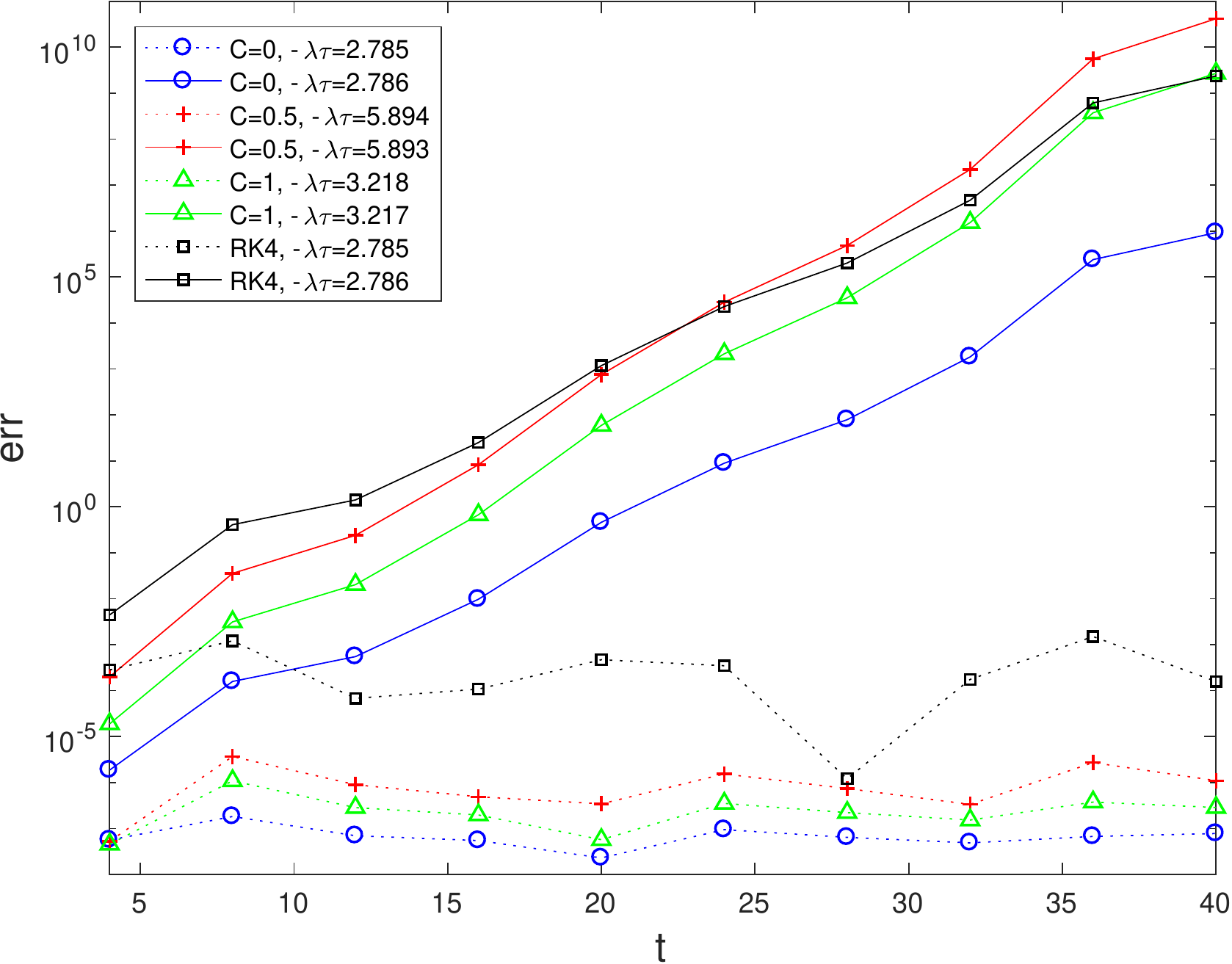}
	\includegraphics[width=0.49\textwidth]{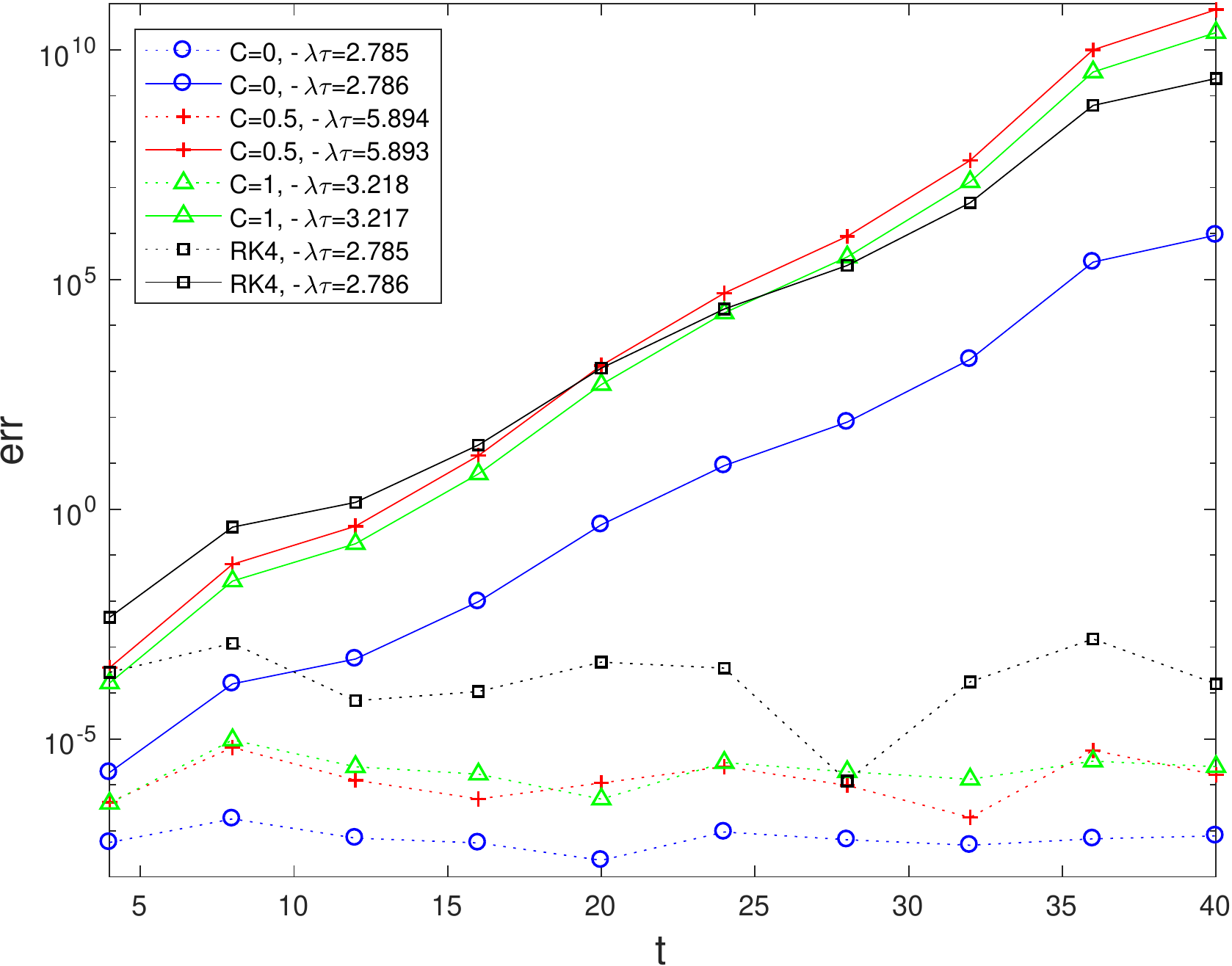}
	\caption{Example \ref{ex002}: Relative errors obtained by the general two-stage fourth-order time discretizations and {\tt RK4} with different step-sizes $\tau$.
		Left: $\alpha = \frac{1}{3}+\frac{C}{60} (\lambda \tau)^3,  ~ \beta = \frac{2}{3}$;
		right: $\alpha = \frac{2}{3},  ~ \beta = \frac{1}{3}+\frac{C}{60} (\lambda \tau)^3$.}
%	\caption{Example \ref{ex002}: errors over time obtained by our methods and {\tt RK4} with different step-sizes $\tau$.
%	Left: $\alpha = \frac{1}{3}+\frac{C}{60} (L_u \tau)^3,  ~ \beta = \frac{2}{3}$;
%		right: $\alpha = \frac{2}{3},  ~ \beta = \frac{1}{3}+\frac{C}{60} (L_u \tau)^3$.}
	\label{fig:ex002}
\end{figure}

\begin{example} \label{ex002}
Solve the initial value problem
\begin{equation*}% \label{eq:ex002}
u'(t) =L(t,u)=\lam (u - \cos t) - \sin t,~ t\geq 0, ~ \lambda=-2100; \quad u(0) = 1,
\end{equation*}
whose  exact solution is $u(t) = \cos t$. In this case,  $L_u = \lam$ is constant so that the step-size
 can  also be taken as a constant following the interval of the absolute stability in Remark \ref{IntervalofRe}.
%Tables \ref{table:ex002-error-0} $-$  \ref{table:ex002-error-2}  list the step-sizes and relative errors obtained by the method \eqref{eq:2stage4ordernew} and {\tt RK4}  with different step-sizes $\tau$.
Figure \ref{fig:ex002} displays the relative errors.
Those results show that
 the errors tend to infinity as time increases if the step-size $\tau$ is chosen as the smallest such that $\lambda \tau\notin I(C)$,
but if the step-size $\tau$ is taken as the biggest such that $\lambda \tau\in I(C)$, then both the present method \eqref{eq:2stage4ordernew} and {\tt RK4} are stable and the errors of \eqref{eq:2stage4ordernew}  are smaller than those of {\tt RK4}. Moreover,
 the biggest step-size for the stability of \eqref{eq:2stage4ordernew} with $C = 0.5$ is almost twice those of \eqref{eq:2stage4ordernew} with $C = 0$ or $1$ and {\tt RK4}.

\end{example}

\begin{figure}[htpb]
	\centering
	\includegraphics[width=0.49\textwidth]{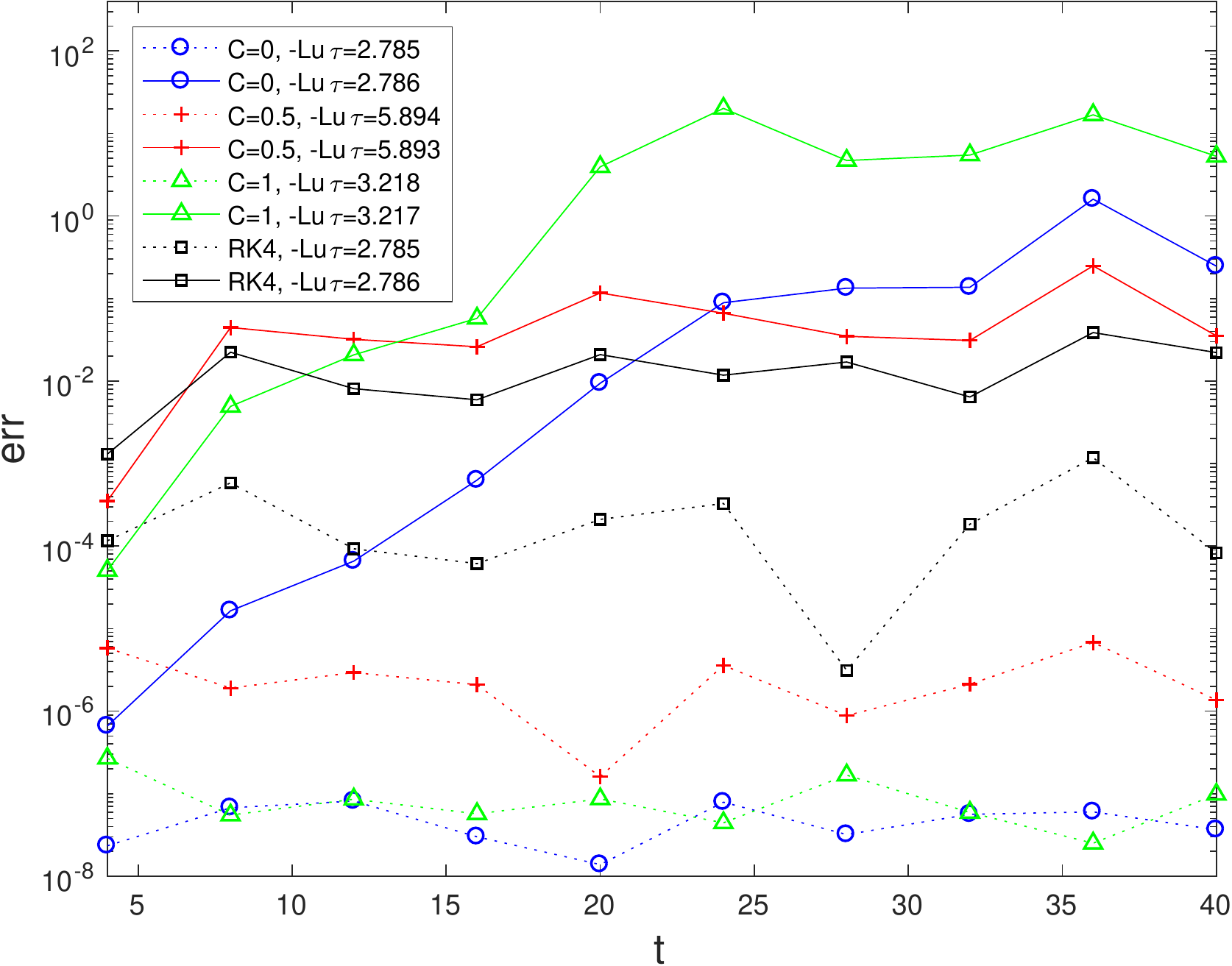}
	\includegraphics[width=0.49\textwidth]{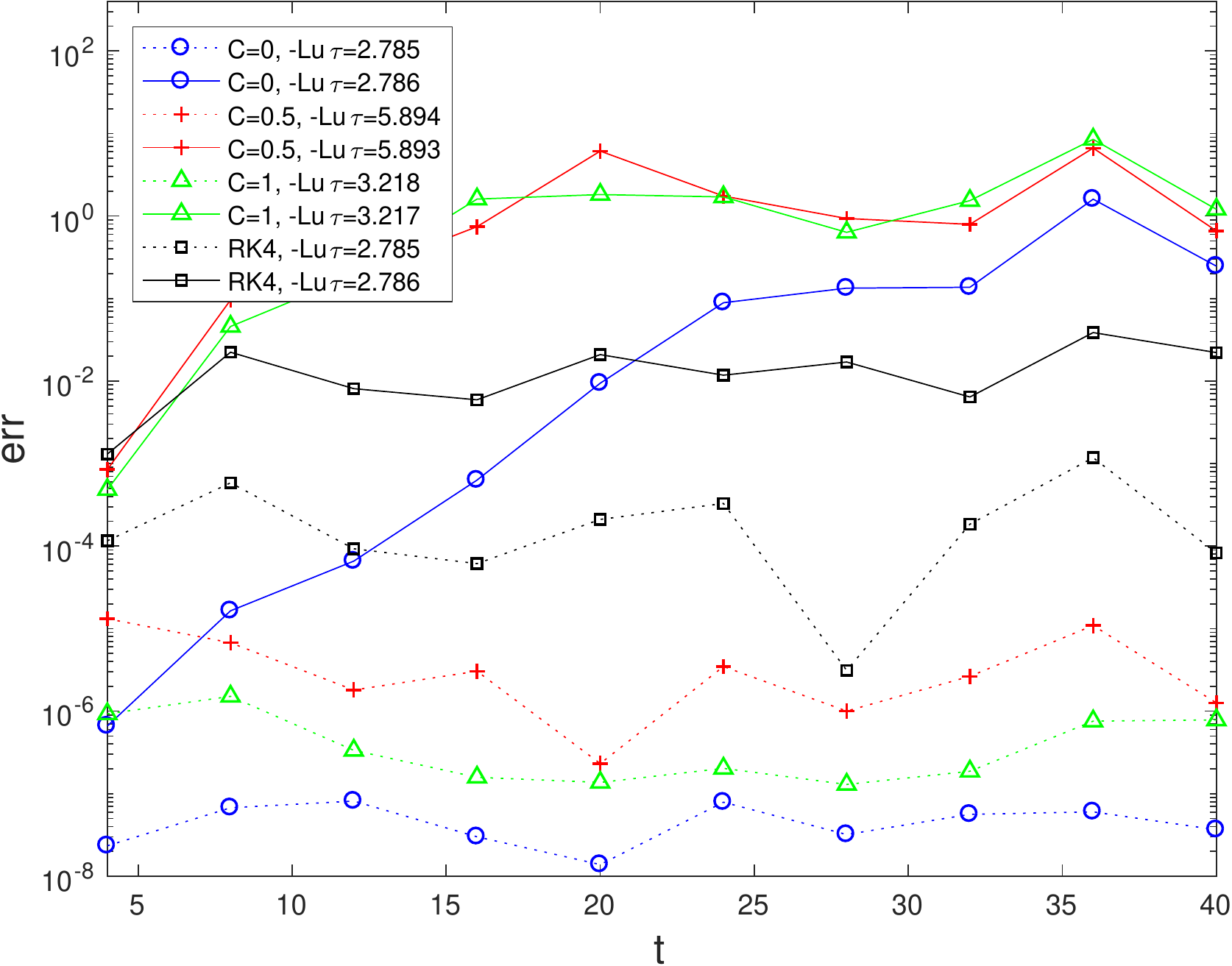}
	\caption{Example \ref{ex002new}: Relative errors obtained by the general two-stage fourth-order time discretizations and {\tt RK4} with different step-sizes $\tau$.
		Left: $\alpha = \frac{1}{3}+\frac{C}{60} (L_u \tau)^3,  ~ \beta = \frac{2}{3}$;
		right: $\alpha = \frac{2}{3},  ~ \beta = \frac{1}{3}+\frac{C}{60} (L_u \tau)^3$.}
	\label{fig:ex002new}
\end{figure}

\begin{example} \label{ex002new}
	Consider the initial value problem of a  nonlinear differential equation
	\begin{equation}\label{eq:ex002new}
	u'(t) = \mu_1(u-\cos(t)) + \mu_2(u^2-\cos^2(t)) - \sin(t), ~ t\geq 0;  \quad u(0) = 1,
	\end{equation}
	 whose   exact solution $u(t) = \cos(t)$.
 In this case, $L_u = \mu_1 + 2\mu_2 u$ is not a constant so that
	the biggest  step-size for the stability is no longer constant.
	
%	Tables \ref{table:ex002new-error-0} -  \ref{table:ex002new-error-2} list the steps and errors obtained by our methods and RK4 with different step-sizes $\tau$. Moreover,
Our calculations take $\mu_1=-2100$ and $\mu_2=10$.
Figure \ref{fig:ex002new} plots the relative errors by our methods and {\tt RK4} with different  $\tau$.
The results show that  the errors of the present method and {\tt RK4} grow over time $t$ if  the step-size $\tau$ is chosen as the smallest such that $L_u \tau\notin I(C)$, but
	if $\tau$ is taken as the biggest such that $L_u \tau\in I(C)$, those time discretizations are  stable and the errors of the proposed methods are smaller than those of {\tt RK4}.
\end{example}

\subsection{System case}
This subsection will solve several system of ordinary differential equations by using the two-stage high-order methods \eqref{eq:2stage4ordernew} with  $\vec \alpha = \frac{1}{3} \vec I+\frac{C\tau^3}{60} (\nabla_{\vec u} \vec L)^3, \beta = \frac{2}{3}$, and $C = 0, 0.5, 1$.

\begin{figure}[htpb]
	\centering
	\includegraphics[width=0.49\textwidth]{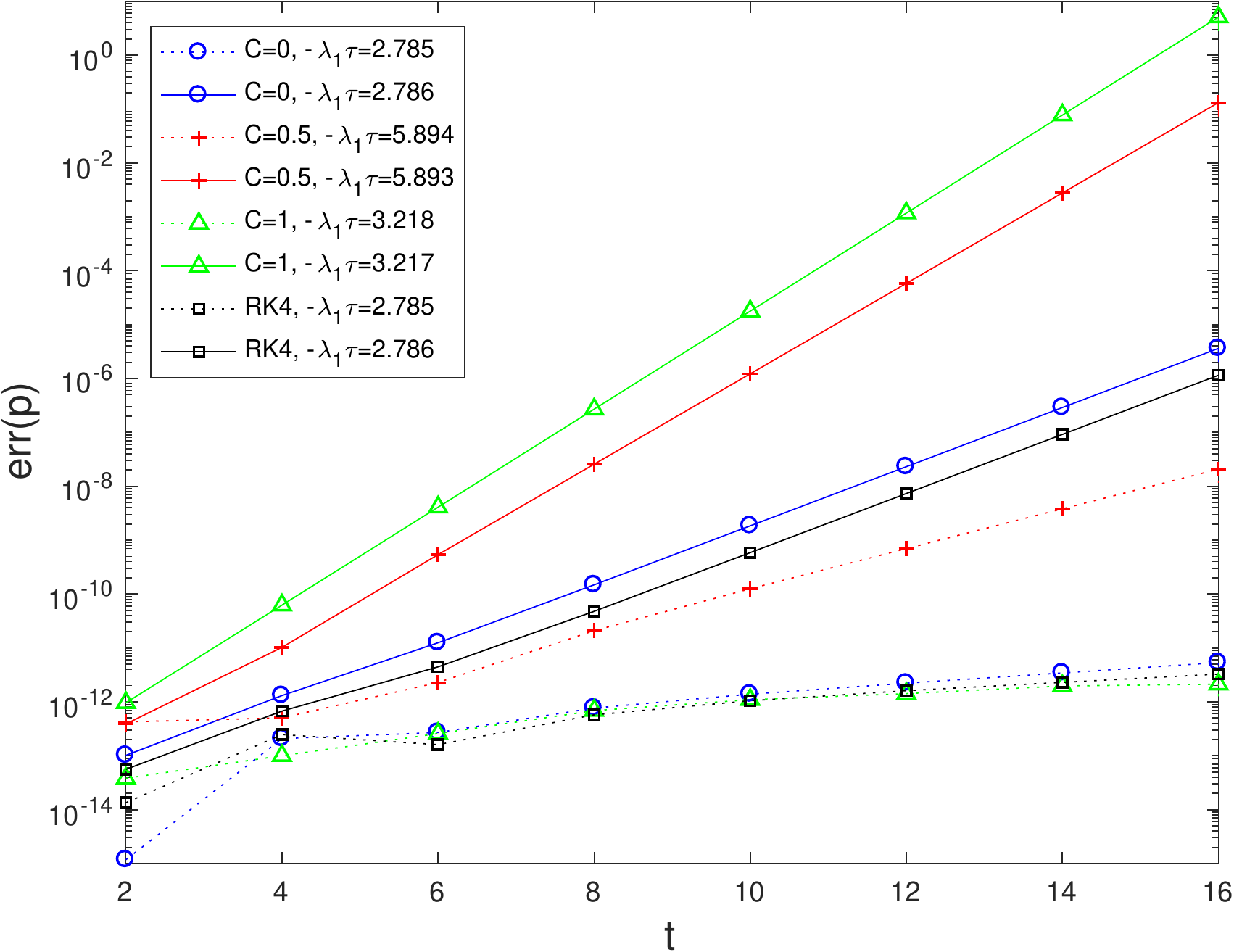}
	\includegraphics[width=0.49\textwidth]{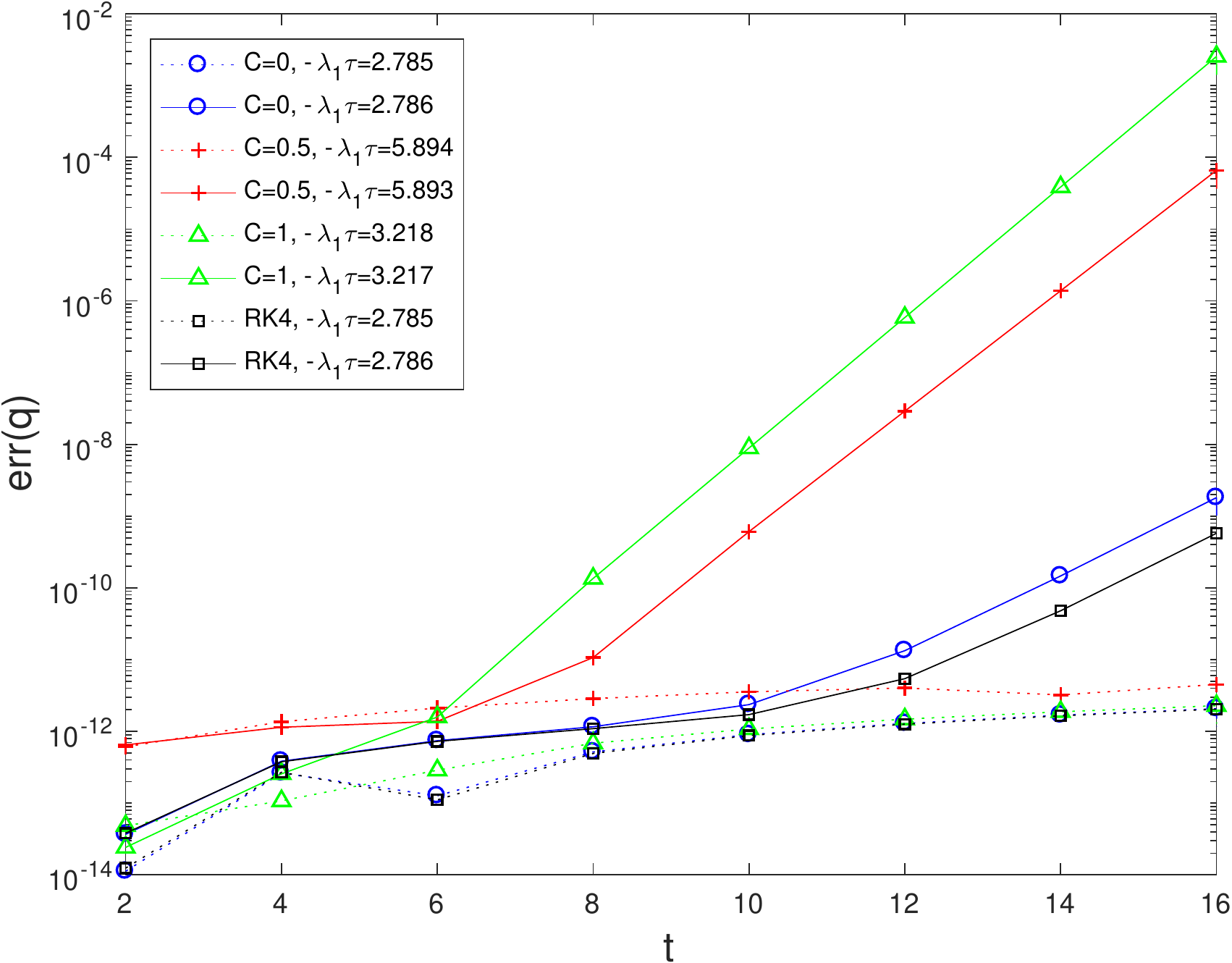}
	\caption{Example \ref{ex006}: Relative errors obtained by the general two-stage fourth-order time discretizations and {\tt RK4} with different step-sizes $\tau$.
		Left: $err(p)$;
		right: $err(q)$.}
	\label{fig:ex006}
\end{figure}

\begin{example} \label{ex006}
This example considers the second-order ODE
$mq''(t) = -kq(t) - cq'(t)$,
describing the motion of a spring oscillator.
By introducing $p(t) = m q'(t)$, it can be converted
to the following system
\begin{equation*} %\label{eq:ex004}
\vec u' = \vec L(\vec u), \quad
\vec u = \begin{pmatrix}
p\\
q
\end{pmatrix},
\quad
\vec L(\vec u) = \begin{pmatrix}
-\frac{c}{m} & -k\\
\frac{1}{m} & 0
\end{pmatrix}
\vec u.
\end{equation*}
Our calculations take the parameters as $m = 1, c = 1001, k = 1000$ and the initial data as $p(0)=-1, q(0)=1$.
Some manipulations show the eigenvalues $\lam_{1} = -1000, \lam_{2} = -1$ of $\vec L_{\vec u}(\vec u)$ and the exact solution $(p, q) = e^{-t} (-1,1) $.

%Table \ref{table:ex006-error-0} -  \ref{table:ex006-error-rk} list the steps and errors at $t = 2,\,4,\,\cdots,\,16$ obtained by our methods and {\tt RK4} with different step-sizes $\tau$.
Figure \ref{fig:ex006} plots the relative errors by our methods and {\tt RK4} with different step-sizes $\tau$.
The results show that the errors of our methods and {\tt RK4} increase to infinity over time if  $\tau$ is chosen as the smallest outside the interval of absolute stability,
but if $\tau$ is chosen as the biggest available step-sizes inside the interval of absolute stability, our methods with $C = 0, 1$ and {\tt RK4} are stable.
Nevertheless, the errors of $p$ obtained by our method with $C = 0.5$ grow over time. %This may be caused by the machine errors.
Table \ref{table:ex006-error-ref} lists the steps and errors obtained by our methods with $C = 0.5,1$ and  $-\lam_{1}\tau = 2.785$.
It can be seen that our methods  are stable now, and the errors are the same size as those with $C = 0$.
\end{example}

\begin{table}[htpb]
	\centering
	\caption{ Example \ref{ex006}: The time steps and errors at $t = 2,\,4,\,\cdots,\,16$ obtained by the general two-stage fourth-order time discretization with $C = 0.5$ or $1$ and $-\lambda_1\tau = 2.785$. }
	\label{table:ex006-error-ref}
	\begin{tabular}{|cc|cc||cc|}
		\hline
		\multirow{2}{*}{t} &  \multirow{2}{*}{{\small step}} &  \multicolumn{2}{c||}{ $C = 0.5$} &   \multicolumn{2}{c|}{ $C = 1$}   \\
		\cline{3-6}
		&  & $err(p)$  & $err(q)$  & $err(p)$  & $err(q)$    \\
		\hline
		2 & 1437 &  2.571e-14 &  2.604e-14 & 5.746e-14 &  5.746e-14 \\
		4 & 2874 &  1.938e-13 &  1.938e-13 & 1.373e-13 &  1.376e-13 \\
		6 & 4311 &  2.325e-13 &  2.325e-13 & 3.114e-13 &  3.113e-13 \\
		8 & 5748 &  6.518e-13 &  6.518e-13 & 7.596e-13 &  7.591e-13 \\
		10 & 7185 &  1.072e-12 &  1.072e-12 & 1.209e-12 &  1.209e-12 \\
		12 & 8622 &  1.492e-12 &  1.492e-12 & 1.655e-12 &  1.655e-12 \\
		14 & 10059 &  1.909e-12 &  1.910e-12 & 2.105e-12 &  2.105e-12 \\
		16 & 11496 &  2.327e-12 &  2.327e-12 & 2.555e-12 &  2.555e-12 \\
		\hline
	\end{tabular}
\end{table}

\begin{figure}[htpb]
\centering
%\subfigure[$x$]{
\includegraphics[width=0.3\textwidth]{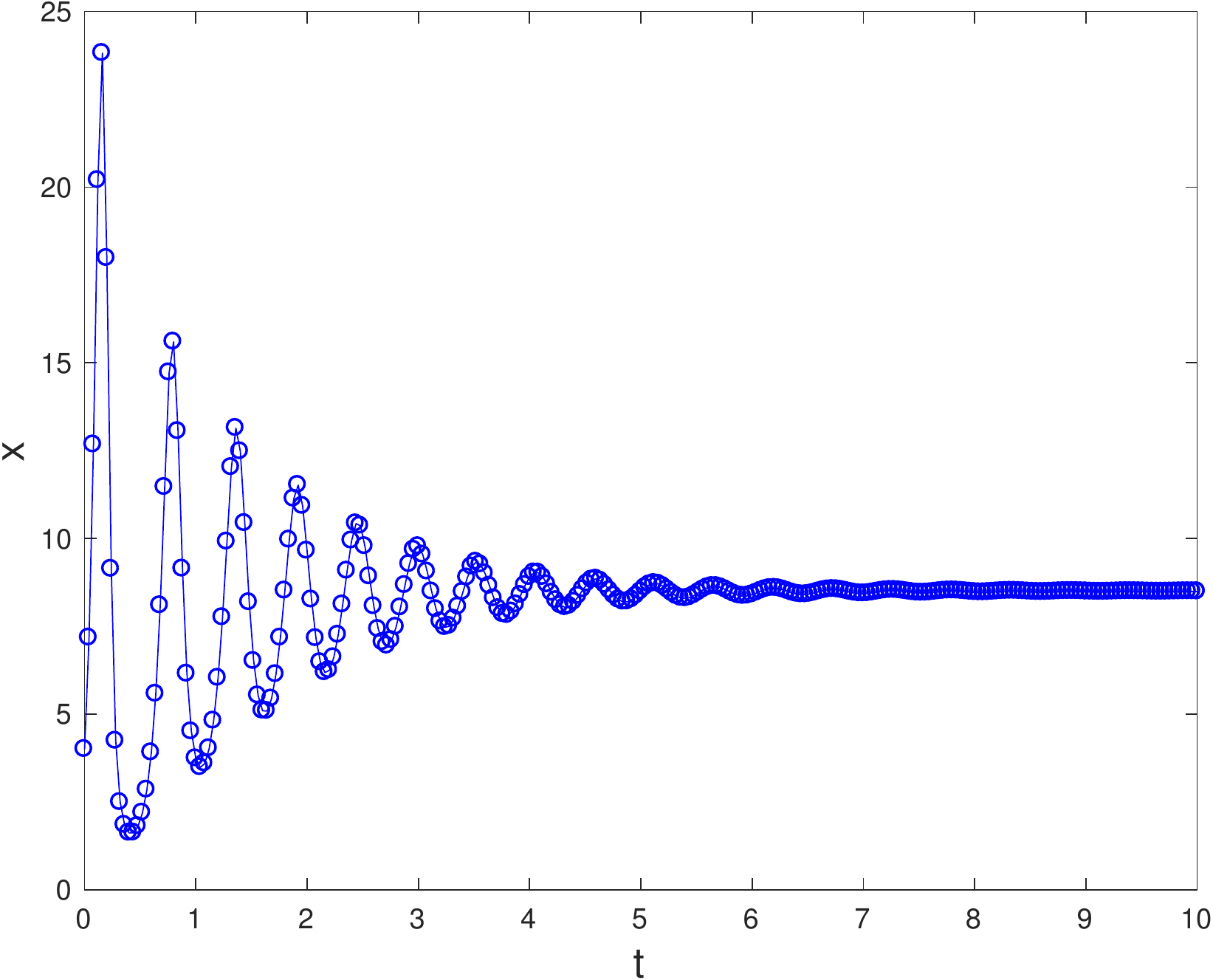}
%}
%\subfigure[$x$]{
\includegraphics[width=0.3\textwidth]{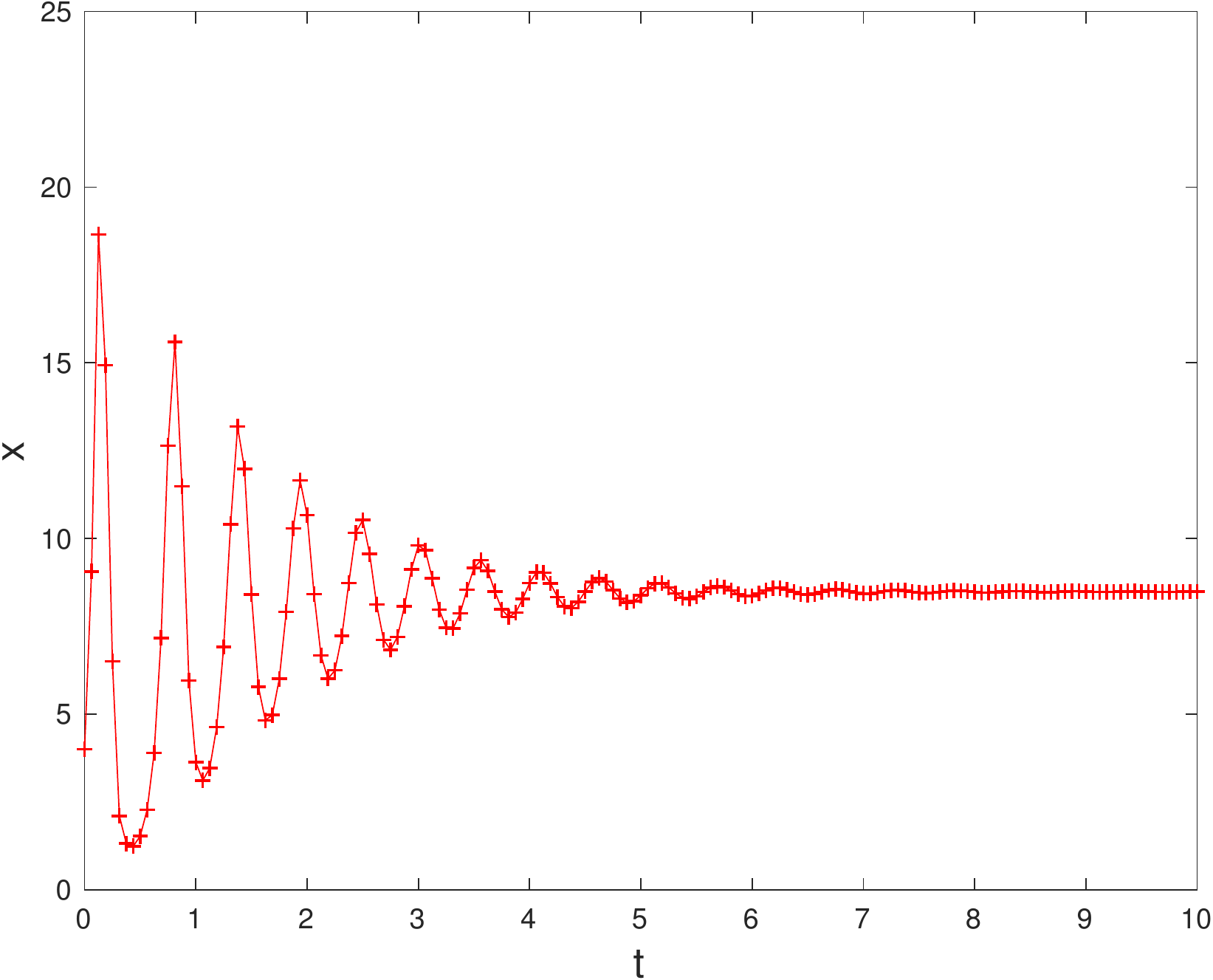}
%}
%\subfigure[$x$]{
\includegraphics[width=0.3\textwidth]{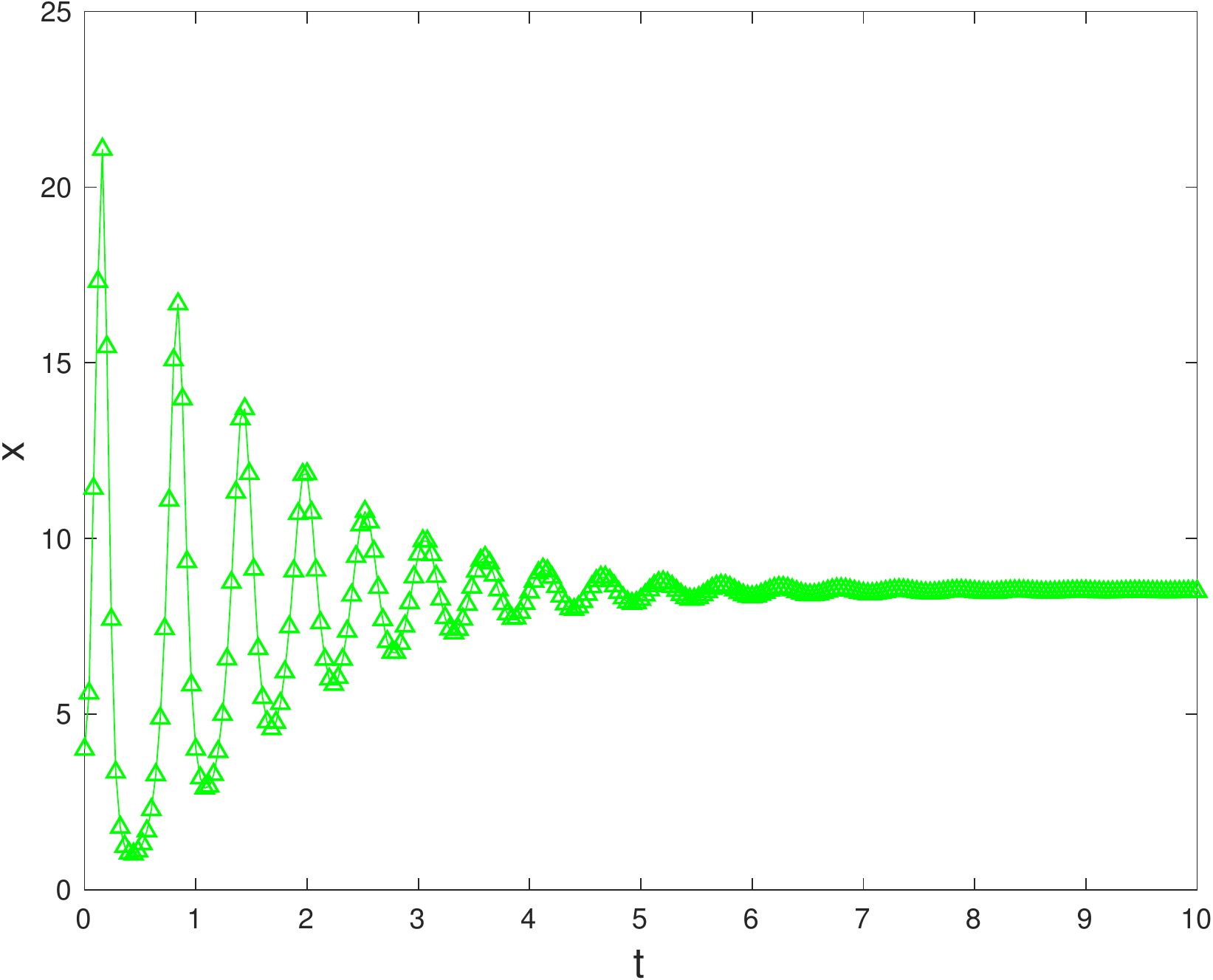}
%}
\\
%\subfigure[$y$]{
\includegraphics[width=0.3\textwidth]{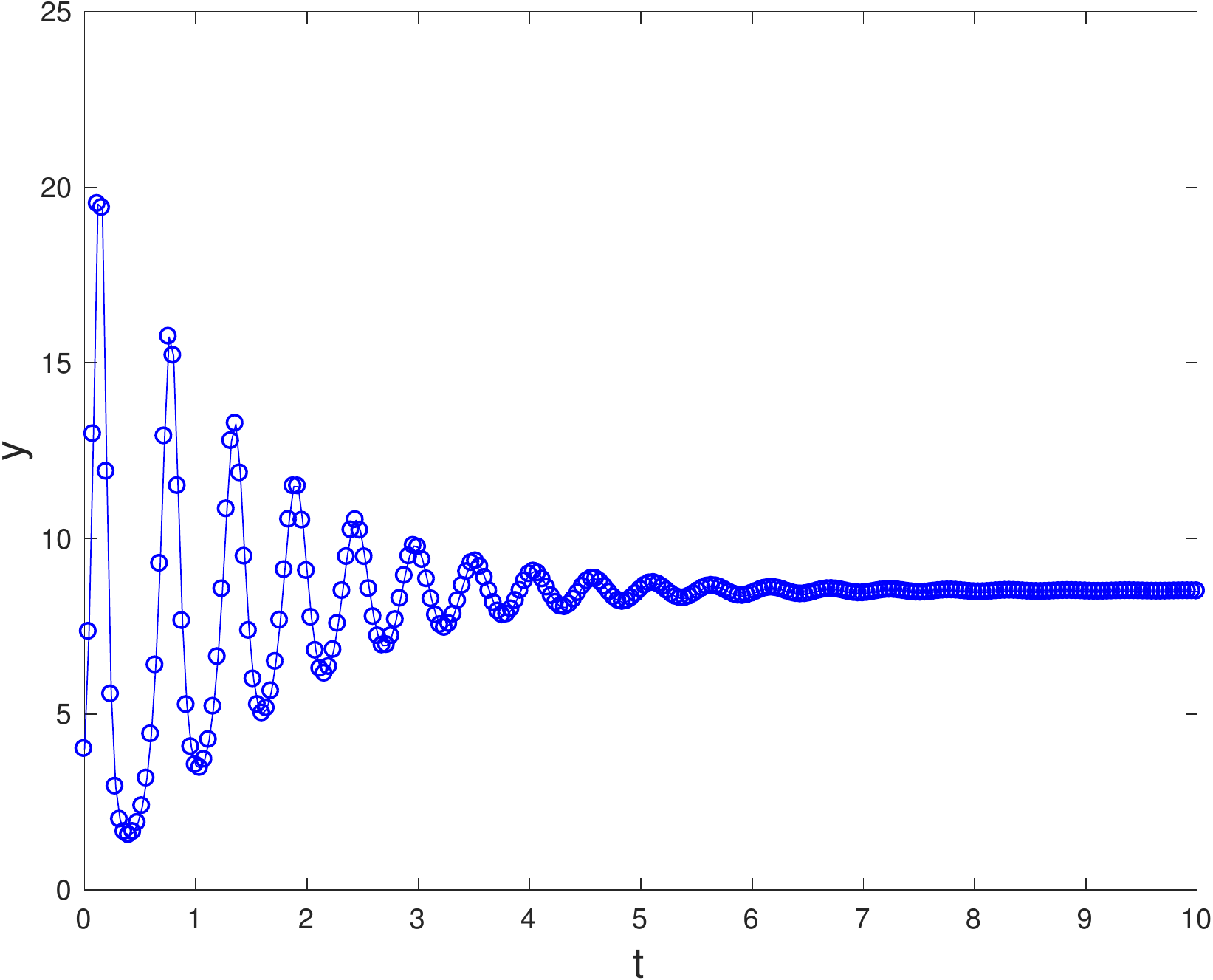}
%}
%\subfigure[$y$]{
\includegraphics[width=0.3\textwidth]{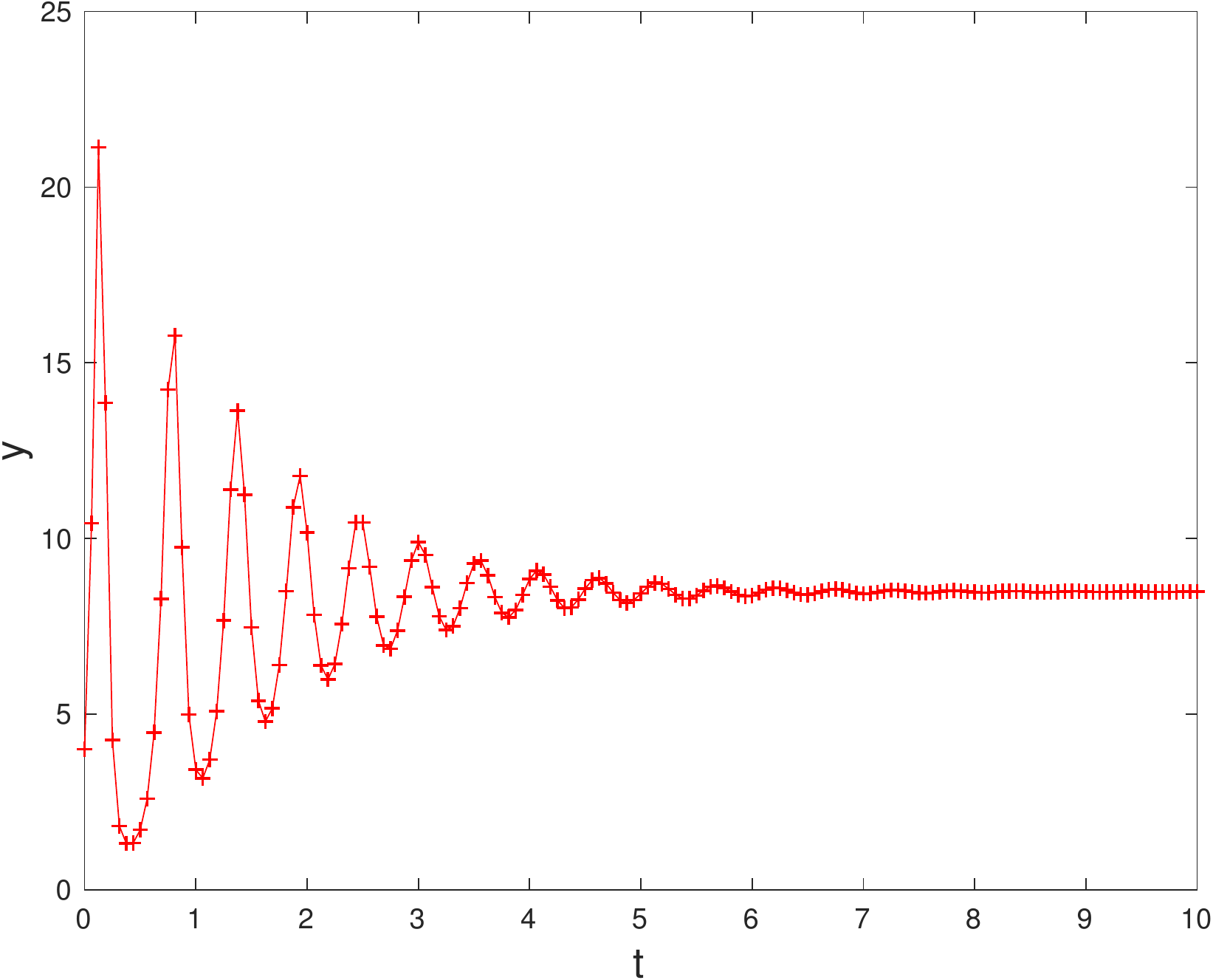}
%}
%\subfigure[$y$]{
\includegraphics[width=0.3\textwidth]{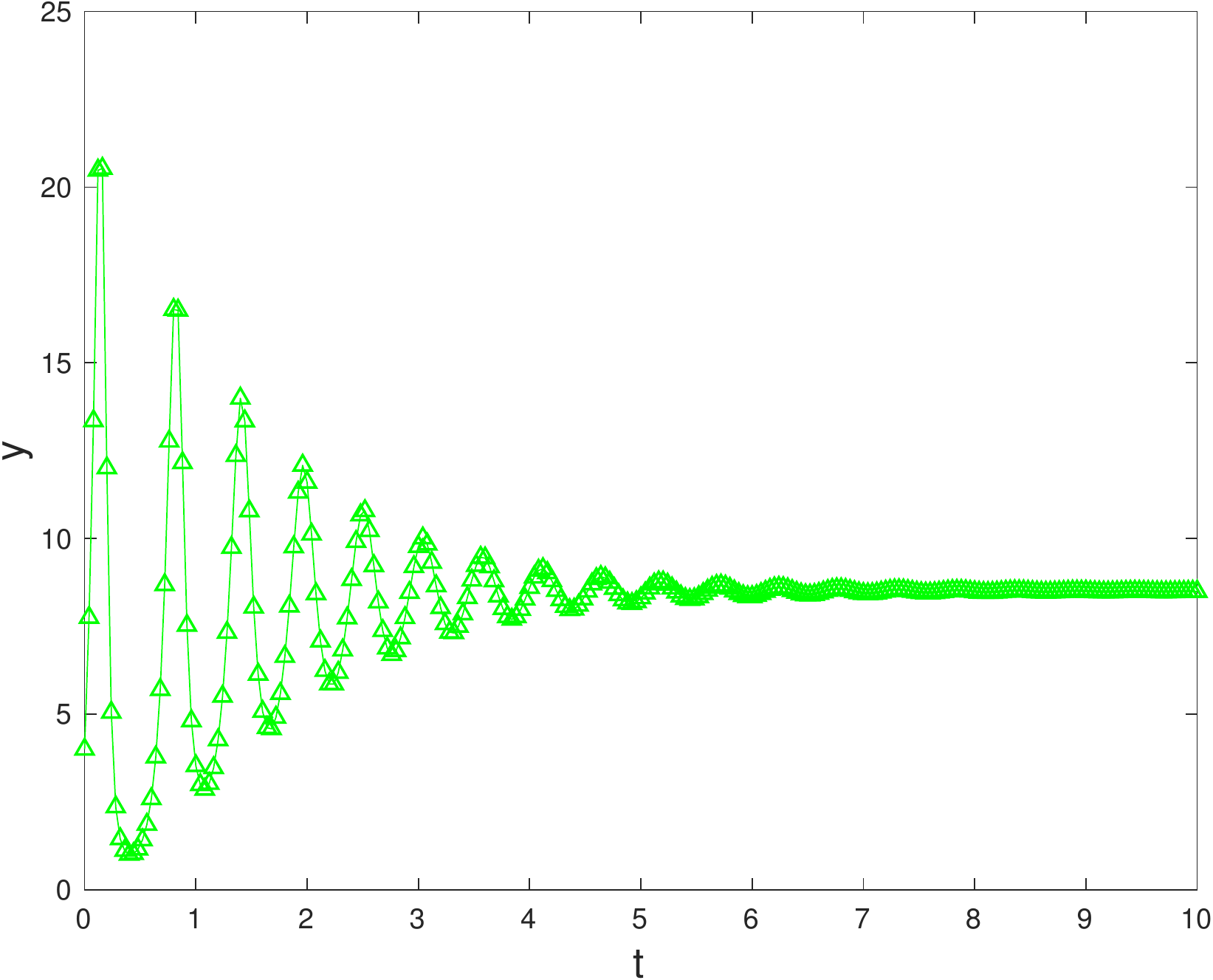}
%}
\\
%\subfigure[$z$]{
\includegraphics[width=0.3\textwidth]{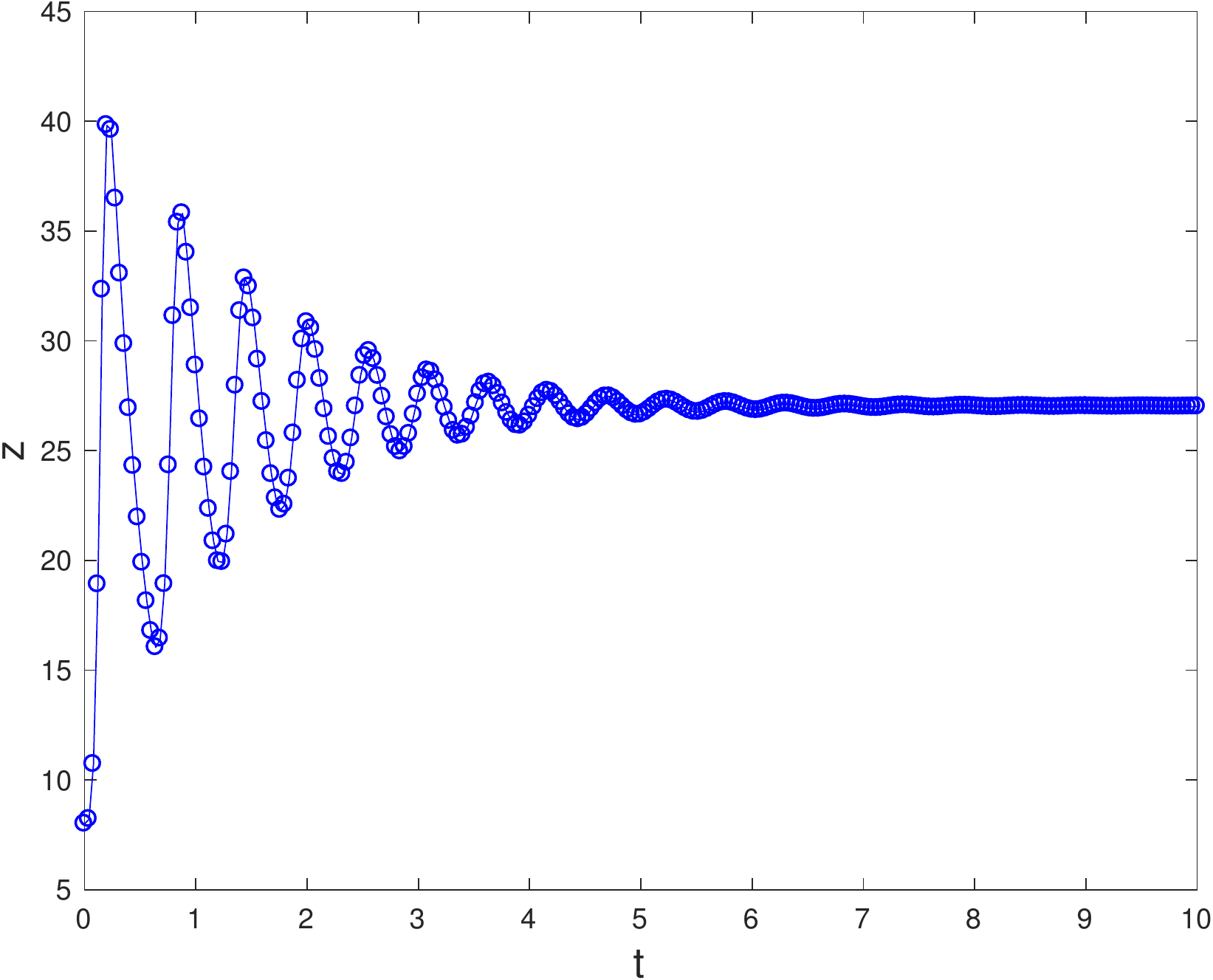}
%}
%\subfigure[$z$]{
\includegraphics[width=0.3\textwidth]{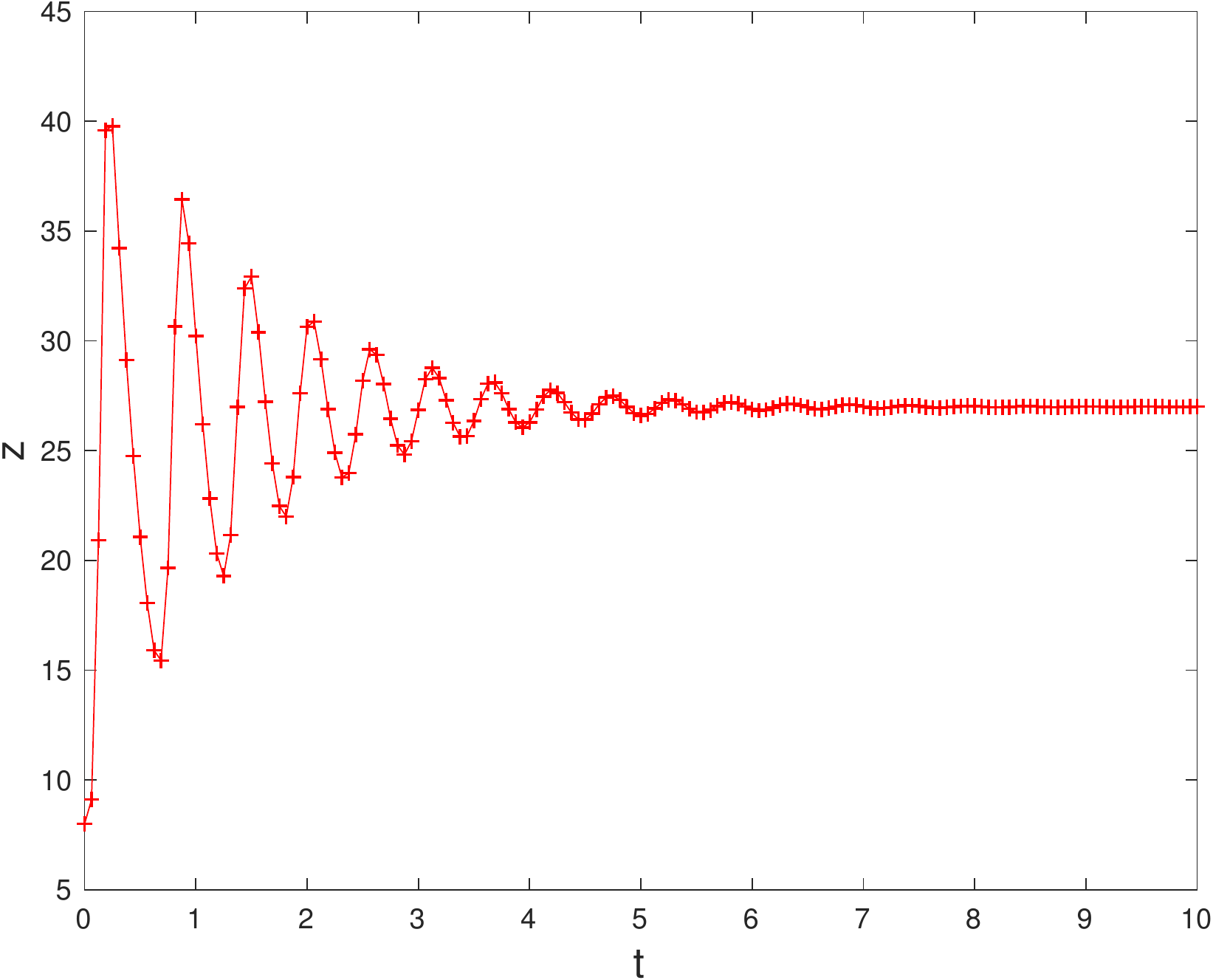}
%}
%\subfigure[$z$]{
\includegraphics[width=0.3\textwidth]{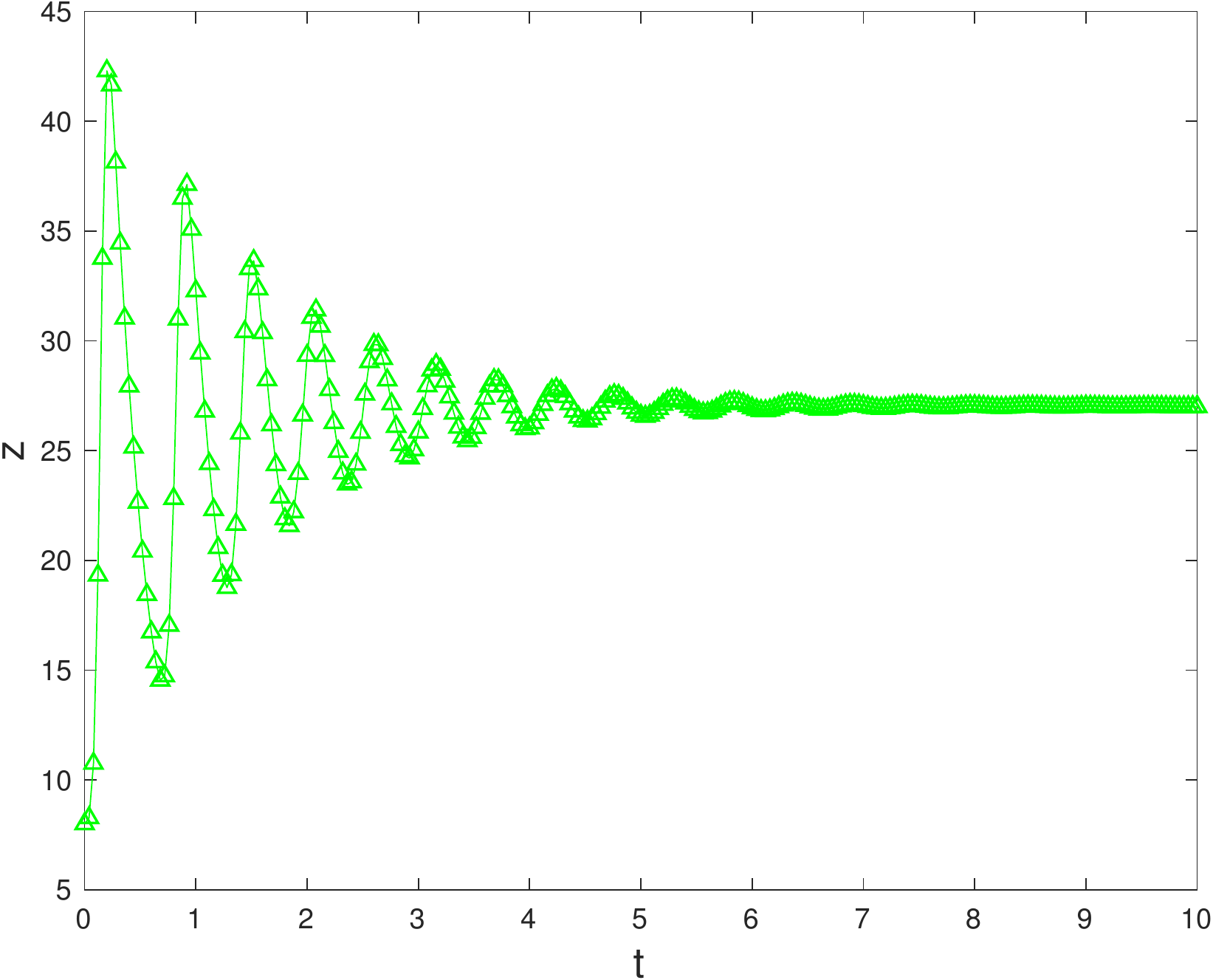}
%}
\\
%\subfigure[Trajectories]{%[Trajectories in phase space $(x,y,z)$]{
\includegraphics[width=0.3\textwidth]{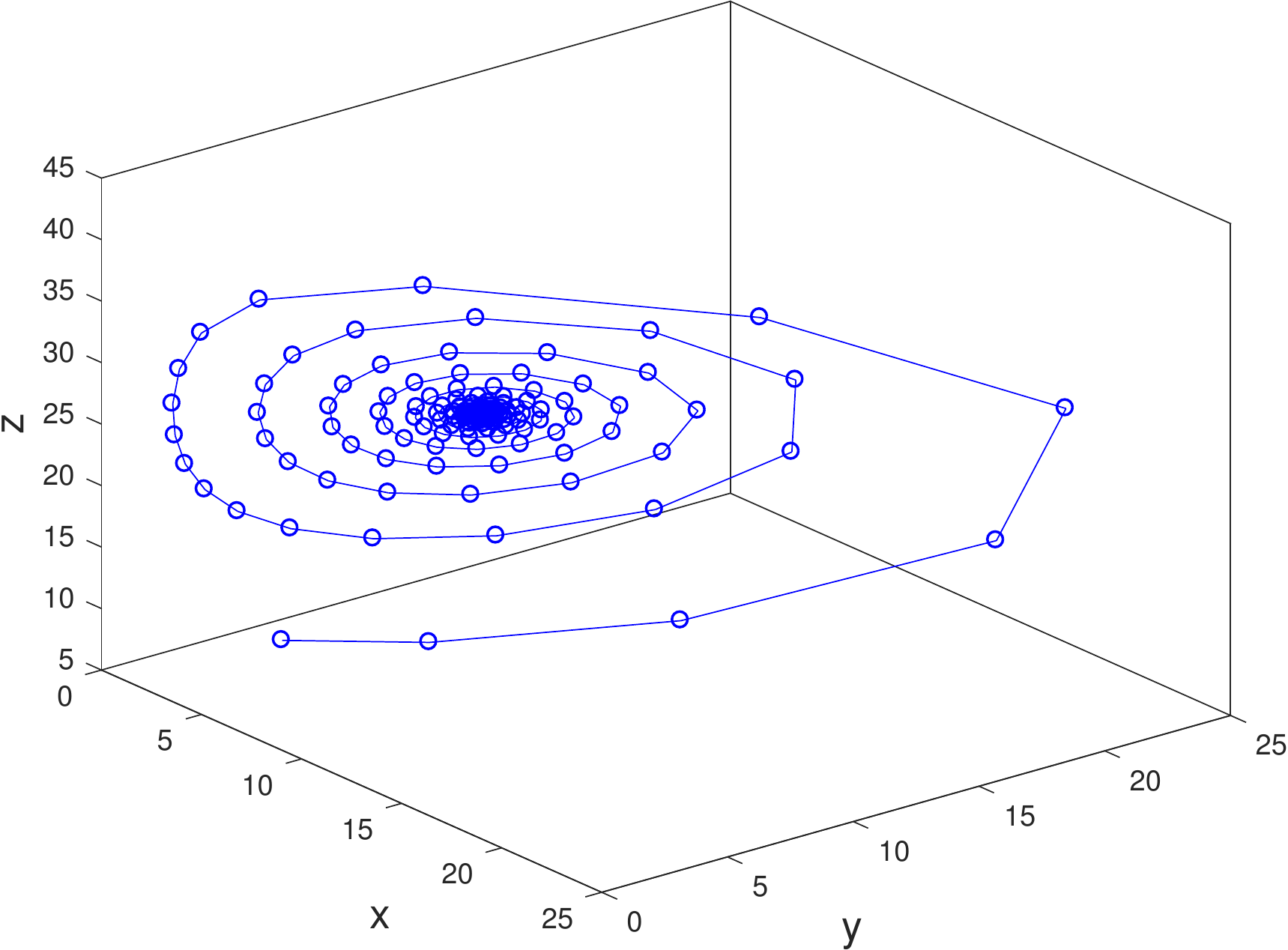}
%}
%\subfigure[Trajectories]{%[Trajectories in phase space $(x,y,z)$]{
\includegraphics[width=0.3\textwidth]{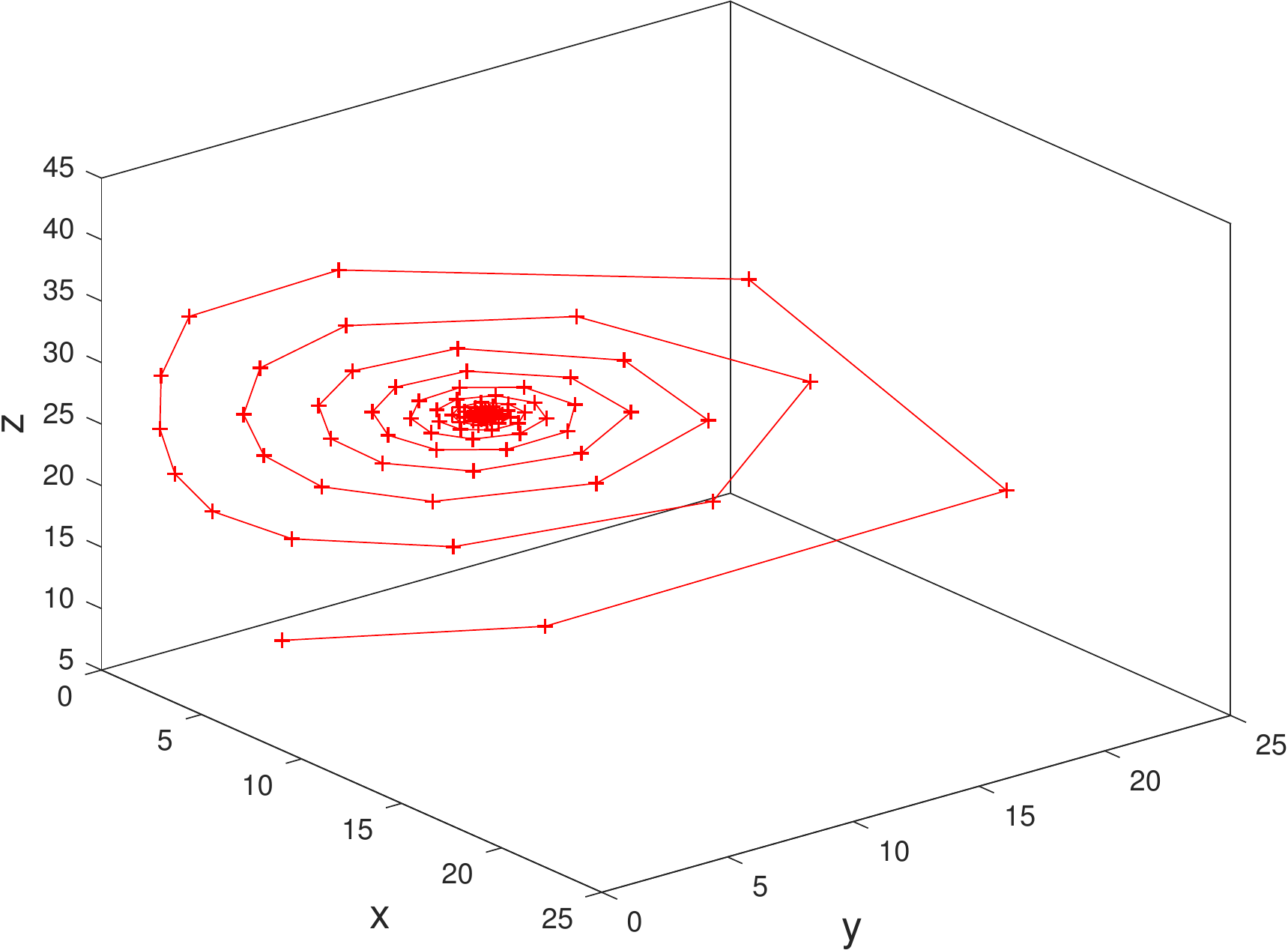}
%}
%\subfigure[Trajectories]{%[Trajectories in phase space $(x,y,z)$]{
\includegraphics[width=0.3\textwidth]{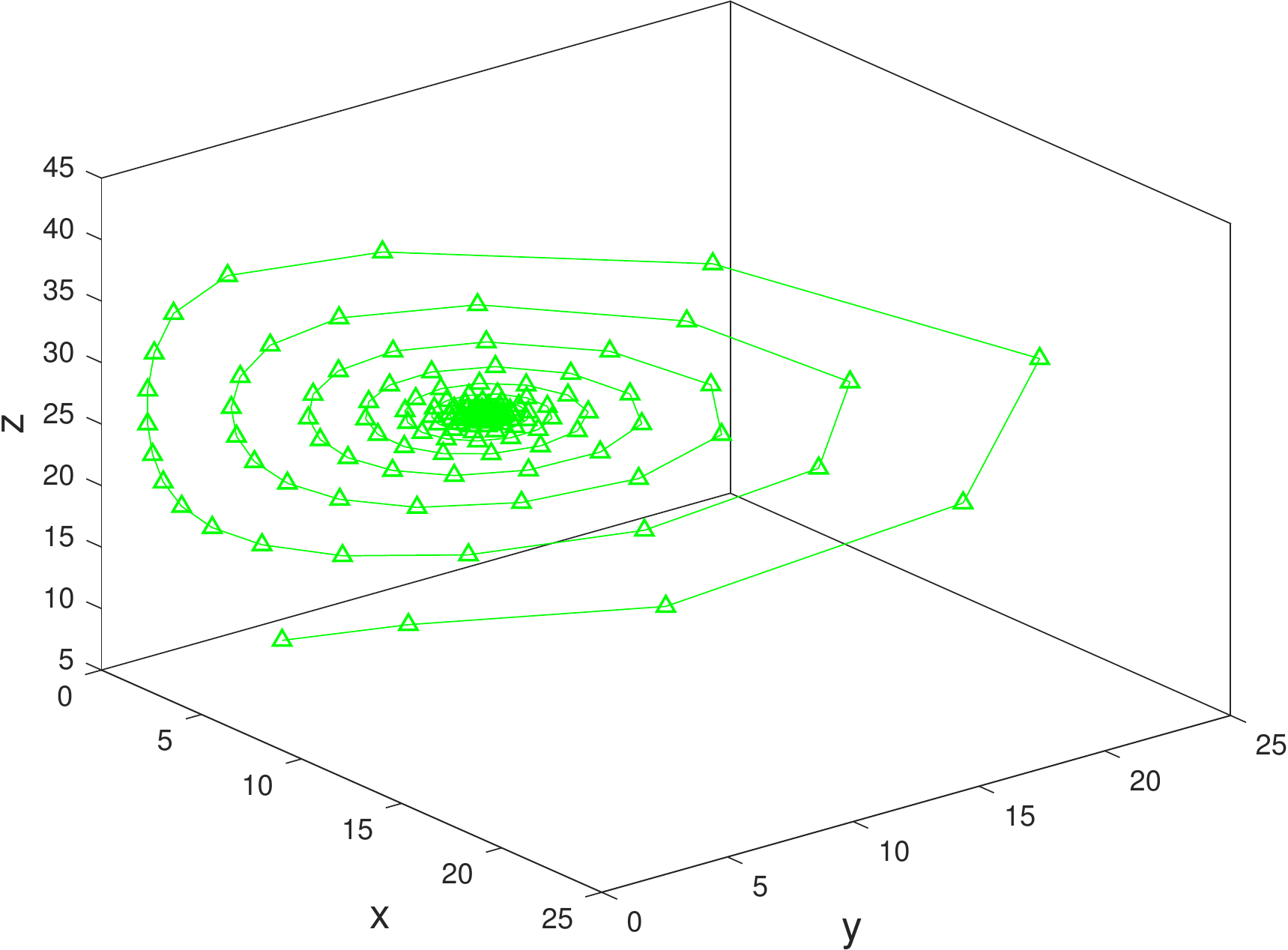}
%}
\caption{Example \ref{ex007}: Solutions (From top to bottom: $x,y,z$ and trajectories) obtained by the general two-stage fourth-order time discretizations. Left: $C = 0, \tau = 0.04$; middle: $C = 0.5, \tau = 0.0625$; right: $C = 1, \tau = 0.04$.}
\label{fig:IndexPa-3-ex007-new}
\end{figure}

\begin{example} \label{ex007}
The last example solves the Lorenz system
\begin{equation*}
\begin{cases}
x'(t) = a(y-x), \\
y'(t) = c x - y - xz, \\
z'(t) = xy - b z,
\end{cases}
\end{equation*}
where $a, b, c$ are constant.
%Lorenz system was first studied by Edward Norton Lorenz in 1963 \cite{Lorenz:1963}.
The behavior depends on the parameters $a,b,c$. %, more details see \cite{Lorenz:1963,Sparrow:1982}.
If we take  $a = 61.8, b = 8/3, c = 28$,  the system has three stationary points: $(x_1,y_1,z_1) = (0,0,0)$, $(x_2,y_2,z_2) = (6\sqrt{2},6\sqrt{2},27)$ and $(x_3,y_3,z_3) = (-6\sqrt{2},-6\sqrt{2},27)$, where the stationary point $(0,0,0)$ is unstable, while the other two are stable.
% then the solutions will go to  $(x_2,y_2,z_2)$ with small enough step-size.

Our calculations take the initial data $x(0) = 4, y(0) = 4, z(0) = 8$.
Tables \ref{table:ex007-new} $-$ \ref{table:ex007-rk-new} list the relative errors
 \begin{equation*}
 err(u) = \frac{|u^{ref}(T) - u_{\tau}(T)|}{|u^{ref}(T)|},\ u=x,y,z, %|u(T) - u_{\tau}(T)|
 \end{equation*}
at $t = 1,2,\cdots,10$ obtained by our methods and {\tt RK4} with different step-sizes $\tau$,
 where $u^{ref}(T)$ is the reference solution obtained by {\tt RK4} with $\tau = 0.001$.
 Figure \ref{fig:IndexPa-3-ex007-new} shows the solutions over time with different $C$ and step-sizes $\tau$.
The results show that the method with $C = 0.5$ permits larger step-size than that with $C = 0$ or $1$.  It is worthy to note that the program of our method with $C = 0$ or $1$ will break up if $\tau = 0.0625$.
\end{example}

%%% new data
\begin{table}[htpb]
	\centering
	\caption{ Example \ref{ex007}: Errors at $t = 1,\,2,\,\cdots,\,10$ obtained by the general two-stage fourth-order time discretization with  $C = 0$.  }
	\label{table:ex007-new}
	\begin{tabular}{|c|ccc|ccc|}
	\hline
	 \multirow{2}{*}{t} &  \multicolumn{3}{c|}{$\tau = 0.04$ } & \multicolumn{3}{c|}{$\tau = 0.01$} \\
	 \cline{2-7}
	 & $err(x)$ & $err(y)$ & $err(z)$ & $err(x)$ & $err(y)$ & $err(z)$ \\
	 \hline
1 &   6.7015e-02 &   2.9769e-03 &   9.9755e-02  &   2.0257e-05  &   1.7648e-05  &   4.4321e-06 \\
 2 &  1.7809e-01 &   2.0776e-01  &  4.4027e-02  &   3.0170e-06  &   5.8543e-06 &    7.1119e-06 \\
 3 &  1.8387e-02 &   5.9763e-03  &  6.1340e-02 &    6.5192e-06  &   4.9609e-06  &   4.4250e-06\\
4   &  4.5944e-02&    3.8950e-02  &  2.0328e-02 &    6.0860e-06  &   5.8296e-06  &   1.0720e-06 \\
 5  & 2.5444e-02 &   2.5753e-02 &   1.3452e-03  &   2.9386e-06  &   3.2706e-06  &   5.3503e-07\\
6  &  7.0759e-03  &  8.7117e-03  &  3.1256e-03 &    4.1393e-07  &   7.6983e-07 &    7.7225e-07\\
 7  &  6.7926e-04 &   3.1301e-04 &   2.2993e-03 &    5.5782e-07 &    3.8339e-07 &    4.4008e-07\\
 8 &   1.8880e-03 &   1.5924e-03 &   8.2682e-04 &    5.3004e-07  &   5.0065e-07 &    1.1080e-07\\
 9  &  1.0623e-03  &  1.0787e-03  &  5.1666e-05 &    2.3573e-07   &  2.6121e-07  &   3.8009e-08\\
 10  & 2.8804e-04 &   3.5932e-04 &   1.3760e-04&     3.0994e-08  &   5.6925e-08  &   5.6218e-08\\
 \hline
	\end{tabular}
\end{table}

\begin{table}[htpb]
	\centering
	\caption{ Example \ref{ex007}: Same as Table \ref {table:ex007-new}
except for $C=0.5$.}
%Errors at $t = 1,\,2,\,\cdots,\,10$ obtained by the general two-stage fourth-order time discretization with  $C=0.5$.  }
	\label{table:ex007-1-new}
	\begin{tabular}{|c|ccc|ccc|}
	\hline
	 \multirow{2}{*}{t} &  \multicolumn{3}{c|}{$\tau = 0.0625$} & \multicolumn{3}{c|}{ $\tau = 0.01$} \\
	 \cline{2-7}
	 & $err(x)$ & $err(y)$ & $err(z)$ & $err(x)$ & $err(y)$ & $err(z)$ \\
	 \hline
1 &  9.3319e-02 &   3.2845e-02 &   5.7565e-02 &    2.0617e-06  &   3.7958e-06 &    4.2273e-06 \\
2&   9.1353e-02&    1.1158e-01  &  3.7513e-02 &    4.8728e-06 &    7.8086e-06 &    6.0057e-06 \\
3&   2.1367e-02  &  9.4735e-03  &  3.4155e-02 &    4.9381e-06  &   3.3219e-06  &   4.0089e-06 \\
4&    2.7067e-02 &   2.4241e-02 &   9.1287e-03 &    4.6001e-06   &  4.3271e-06  &   1.0246e-06 \\
5&   1.2676e-02  &  1.3233e-02  &  2.5175e-04  &   2.2045e-06  &   2.4134e-06  &   2.9472e-07 \\
6&    2.8142e-03 &   3.7420e-03  &  1.8532e-03 &    3.6937e-07  &   6.0478e-07  &   5.0476e-07  \\
7&    6.7871e-04  &  2.0339e-04 &   1.1253e-03  &   3.2916e-07  &   2.1215e-07  &   2.9230e-07 \\
8&    9.6442e-04  &  8.4702e-04&    3.4594e-04 &    3.3149e-07  &   3.0952e-07 &    7.7360e-08 \\
9&   4.7010e-04 &   4.8793e-04  &   1.1421e-06  &   1.5024e-07  &   1.6457e-07   &  1.9940e-08 \\
10&  1.0853e-04 &   1.4200e-04 &   6.6884e-05  &   2.2575e-08  &   3.8093e-08 &    3.3305e-08 \\
 \hline
	\end{tabular}
\end{table}

\begin{table}[htpb]
	\centering
	\caption{ Example \ref{ex007}: Same as Table \ref {table:ex007-new}
except for $C=1$.}
	\label{table:ex007-2-new-1}
	\begin{tabular}{|c|ccc|ccc|}
	\hline
	 \multirow{2}{*}{t} &  \multicolumn{3}{c|}{$\tau = 0.04$} & \multicolumn{3}{c|}{$\tau = 0.01$} \\
	 \cline{2-7}
	 & $err(x)$ & $err(y)$ & $err(z)$ & $err(x)$ & $err(y)$ & $err(z)$ \\
	 \hline
1&    1.5361e-03 &  1.9805e-03  &  6.4616e-03  &  1.6156e-05  &  1.0065e-05 &   4.0029e-06\\
2&    8.5945e-03 &   1.3262e-02  &  7.0901e-03 &   6.7043e-06 &   9.7256e-06 &   4.8799e-06\\
3&    4.9792e-03 &   3.1234e-03 &   4.8734e-03 &   3.3433e-06 &   1.6740e-06 &   3.5795e-06\\
4&    4.4063e-03  &  4.1019e-03 &   1.1500e-03  &  3.1021e-06 &   2.8133e-06  &  9.7422e-07\\
5&    1.9307e-03 &   2.0774e-03 &   1.7791e-04  &  1.4654e-06  &  1.5507e-06 &   5.3838e-08\\
 6&   3.4475e-04 &   5.1191e-04 &   3.5437e-04  &  3.2398e-07  &  4.3844e-07 &   2.3630e-07\\
7&    1.8519e-04  &  1.0561e-04  &  1.9561e-04  &  9.9956e-08  &  4.0592e-08 &   1.4400e-07\\
8&    1.9198e-04 &   1.7594e-04  &  5.2216e-05  &  1.3241e-07 &   1.1792e-07 &   4.3785e-08\\
9&    8.5114e-05  &  9.1504e-05  &  7.3770e-06 &   6.4516e-08 &   6.7691e-08  &  1.8473e-09\\
10&   1.4865e-05 &   2.2347e-05  &  1.5720e-05 &   1.4119e-08 &   1.9205e-08 &   1.0349e-08\\
 \hline
	\end{tabular}
\end{table}

\begin{table}[htpb]
	\centering
	\caption{ Example \ref{ex007}: Same as Table \ref {table:ex007-new}
except for {\tt RK4}.  }
	\label{table:ex007-rk-new}
	\begin{tabular}{|c|ccc|ccc|}
	\hline
	 \multirow{2}{*}{t} &  \multicolumn{3}{c|}{ $\tau = 0.04$} & \multicolumn{3}{c|}{$\tau = 0.01$} \\
	 \cline{2-7}
	& $err(x)$ & $err(y)$ & $err(z)$ & $err(x)$ & $err(y)$ & $err(z)$ \\
	 \hline
1&     4.2184e-02  & 2.3244e-02  & 2.1487e-02  &   4.0999e-05  &   2.2965e-05 &    8.8756e-06\\
2&     2.1815e-02 &  3.3483e-02  & 2.3926e-02 &    1.0701e-05  &   1.6094e-05&     1.0973e-05\\
3&     1.7573e-02 &  1.3117e-02  & 1.3310e-02 &    8.9941e-06 &    6.7554e-06 &    7.6197e-06\\
4&     1.3504e-02 &  1.2992e-02 &  2.3494e-03 &     8.8726e-06 &    8.4220e-06 &    1.7920e-06\\
5&     5.1771e-03  & 5.8316e-03  & 1.0917e-03  &   4.1519e-06 &    4.5617e-06 &    6.5426e-07\\
6&     4.3010e-04 &  9.9634e-04 &  1.2432e-03 &    6.1669e-07 &    1.0786e-06 &    9.9454e-07\\
7&     8.7431e-04 &  6.4763e-04  & 5.8343e-04  &   6.7845e-07  &   4.5406e-07 &    5.6068e-07\\
8&     6.6585e-04 &  6.4333e-04 &  1.0654e-04 &    6.4966e-07 &    6.1090e-07 &    1.4192e-07\\
9&     2.4192e-04 &  2.7694e-04  & 5.9139e-05  &   2.8740e-07  &   3.1716e-07  &   4.3333e-08\\
10&    1.1922e-05  & 3.9662e-05  & 6.1925e-05 &    3.9208e-08 &    6.9809e-08  &   6.6129e-08\\
 \hline
	\end{tabular}
\end{table}

\section{Conclusion} \label{sec:conclusion}
By introducing variable weights, this paper proposed  a class of more general explicit one-step two-stage time discretizations, which are different from the existing methods, such as the Euler methods, Runge-Kutta methods, and multistage multiderivative methods etc.
Their absolute stability, the stability interval, and the intersection between the imaginary axis and the absolute stability region were carefully studied. The results showed that the new two-stage time
discretizations could be fourth-order accurate conditionally, the absolute stability region of the proposed methods with some special choices of the variable weights could be larger than
that of the classical explicit fourth- or fifth-order Runge-Kutta method, and the interval of
absolute stability can be almost twice as much as the latter.
Several numerical experiments
 were carried out to demonstrate the performance and accuracy as well as the stability of the
proposed methods. It is interesting to apply the present
 time discretization to solving the time-dependent partial differential equations.

%Under certain conditions those two-stage time discretizations were fourth-order accurate, and their absolute stability regions would change with the choices of the parameters  in \eqref{eq:2stage4ordernew}.
%More importantly, for some special choices, for example $\alpha + \beta = 1 + \frac{C}{60} (L_u \tau)^3$ and $C = 0.5$, the absolute stability region and interval
%were larger than the classical fourth- or fifth-order Runge-Kutta method.
%%Specially, the absolute stability region could be enlarged after some careful choice -- two-stage fourth-order methods \eqref{eq:2stage4ordernew} with
%%In other words, one may have the  explicit two-stage fourth-order time discretizations with a more larger
%%time step-size.
%Several  numerical experiments showed that our proposed methods with $\alpha + \beta = 1 + \frac{C}{60} (L_u \tau)^3$ and $C = 0, 0.5, 1$ attained the theoretical order.
%Moreover, our method with $\alpha = \frac{1}{3}+\frac{C}{60} (L_u \tau)^3,~\beta = \frac{2}{3}$ and $C = 0.5$ permitted larger time-step than that with  $C = 0, 1$ and method {\tt RK4}.

%% Acknowledgments %%%%%%%%
\section*{Acknowledgments}
The authors were partially supported by the Special Project on High-performance Computing under the
National Key R\&D Program (No. 2016YFB0200603), Science Challenge Project (No. TZ2016002), the Sino-German Cooperation Group Project (No. GZ 1465),
and the National Natural Science Foundation of China (No. 11421101).
%and High-performance Computing Platform of Peking University.
%http://www.sinogermanscience.org.cn/de/aktuelles/de_2018/201807/t20180726_31020.html

%\section*{\refname}
%\bibliography{journalname,ASRgeneral}
%%% Bibliography  %%%%%%%%%%

\end{document}